\documentclass{article}

\usepackage[utf8]{inputenc}
\usepackage[english]{babel}

\usepackage{graphicx}
\usepackage{amssymb}
\usepackage{amsmath}
\usepackage{tikz-cd}
\usepackage{enumerate}
\usepackage{verbatim}
\usepackage{mathtools}
\usepackage{relsize}

\usepackage{hyperref}

\newtheorem{theorem}{Theorem}[section]
\newtheorem{corollary}[theorem]{Corollary}
\newtheorem{lemma}[theorem]{Lemma}
\newtheorem{proposition}[theorem]{Proposition}
\newtheorem{conjecture}[theorem]{Conjecture}

\newcommand{\Rep}{\text{Rep}}

\begin{document}

\title{Effective and Infinite-Rank Superrigidity in the Context of Representation Stability}
\author{{Nate Harman}}


\maketitle

\vspace{-.5cm}
\begin{abstract}
We discuss certain effective improvements on superrigidity for $SL_n(\mathbb{Z})$ for finite $n>2$. Using these ideas we then use superrigidity to prove a representation stability theorem about pointwise finite dimensional $VIC(\mathbb{Z})$-modules, which itself can be viewed as a superrigidity theorem for $VIC(\mathbb{Z})$ and $GL_\infty(\mathbb{Z})$.

\end{abstract}

\section{Introduction}
Superrigidity for $SL_n(\mathbb{Z})$ with $n>2$, is the statement that any finite dimensional representation of $SL_n(\mathbb{Z})$ (or any finite index subgroup thereof) virtually extends to an algebraic (a.k.a. rational) representation of $SL_n(\mathbb{C})$. Explicitly this means that for any representation $$\phi: SL_n(\mathbb{Z}) \to GL_N(\mathbb{C})$$ there exists a finite index subgroup $\Gamma \subseteq SL_n(\mathbb{Z})$  and an algebraic representation $$ \tilde{\phi}: SL_n(\mathbb{C}) \to GL_N(\mathbb{C})$$ such that $\phi(g) = \tilde{\phi}(g)$ for all $g \in \Gamma$.  This statement was first proved by Bass, Milnor, and Serre (\cite{BMS} Theorem 16.2) although the term ``superrigidity" is usually associated with Margulis, who proved a much more general version for arbitrary irreducible lattices in higher rank semisimple Lie groups (see \cite{marg}).

From a representation theory perspective, superrigidity gives us a fairly complete understanding of the finite dimensional representations of $SL_n(\mathbb{Z})$ in terms of the algebraic representations of $SL_n(\mathbb{C})$ and the representations of the finite quotients $SL_n(\mathbb{Z}/\ell \mathbb{Z})$. In the first part of this paper (Section \ref{fixedn}) we give a largely expository account of this theory, with a particular focus on what the representation theoretic properties of a representation can tell us about which finite index subgroup actually arises in the statement of superrigidity, and discuss ``effective" improvements on superrigidity such as:

\medskip

\noindent \textbf{Corollary 2.8} \emph{ For $n \ge 5$ any representation of $SL_n(\mathbb{Z})$ of dimension less than $2^n-2$ extends (not just virtually) to an algebraic representation of $SL_n(\mathbb{C})$.}

\medskip

Representation stability, broadly speaking, refers to a phenomenon in which certain sequences of representations $V_n$ of sequences of groups $G_n$ appear to ``stabilize" in $n$, despite being representations of different groups.  The most well known part of this theory is the theory of finitely generated $FI$-modules first developed by Church, Ellenberg, and Farb \cite{CEF}.  A finitely generated $FI$-module gives rise to a sequence of vector spaces $V_1,  V_2, V_3, \dots$ with $V_n$ carrying an action of the symmetric group $S_n$, which stabilize in $n$ in an appropriate sense.

$VIC(R)$-modules for $R$ a commutative ring were defined by Putman and Sam in \cite{PS}, they are analogs of $FI$-modules with general linear groups $GL_n(R)$ playing the same role that symmetric groups played for $FI$-modules (see Section \ref{vicsect} for a precise definition). The most well-understood case is when $R = \mathbb{F}_q$ is a finite field (studied in \cite{PS}, \cite{GW}, and \cite{MW}), where we now know many of the same properties as for $FI$-modules.

Much less well understood, but important for many applications, is the case when $R = \mathbb{Z}$ where now we obtain compatible sequences of $GL_n(\mathbb{Z})$-representations. These have arisen in the study of the Torelli group $IA_n$ defined as the kernel of the map natural map $$\text{Aut}(F_n) \to GL_n(\mathbb{Z})$$ coming from identifying the abelianization of $F_n$ with $\mathbb{Z}^n$. $GL_n(\mathbb{Z})$ acts on the homology of $IA_n$ and in fact for fixed $k$ the homology groups $H_k(IA_n)$ form a $VIC(\mathbb{Z})$-module. While this and related constructions are certainly a major motivation for this work, we won't say anything further about these. For more information about these we'll refer to \cite{MPW} and \cite{DP}.

The main objects of study in this paper will be finitely generated $VIC(\mathbb{Z})$-modules over the complex numbers with the additional constraint that all of the representations $V_n$ are finite dimensional (a condition that comes for free for finitely generated $FI$ or $VIC(\mathbb{F}_q)$-modules). Not much is known previously about these modules, but we would like to highlight two known results as they pertain to this work:

\begin{enumerate}

\item Putman and Sam have shown that, without the finite dimensionality constraint, finitely generated $VIC(\mathbb{Z})$-modules can fail to be Noetherian (\cite{PS} Theorem N), meaning they can have submodules which are not finitely generated. This is a problem for many of the spectral sequence arguments that are common in the theory of representation stability. Our Theorem \ref{Noetherianity} says that under this additional finite dimensionality condition Noetherianity does hold.

\item In \cite{Patzt1} Patzt extensively looked at the class of $VIC(\mathbb{C})$-modules for which the representations $V_n$ are pointwise algebraic (actually he works with $VIC(\mathbb{Q})$, but there is no difference for pointwise algebraic modules). He has proved a Noetherianity statement as well as a multiplicity stability statement for such modules. His results will be a special case of the theory developed here, but we will be drawing heavily on his ideas to understand the ``algebraic parts" of $VIC(\mathbb{Z})$-modules.

\end{enumerate}

The main part of this paper, especially Section \ref{vicsect}, will be about the interaction of superrigidity and representation stability.  Superrigidity will let us control the behavior of finitely generated $VIC(\mathbb{Z})$-modules with our main theorem, Theorem \ref{main}, stating that finitely generated, pointwise finite dimensional $VIC(\mathbb{Z})$-modules over the complex numbers are in a sense built out of a stable algebraic part, and finitely generated $VIC(\mathbb{Z}/\ell \mathbb{Z})$-modules.  Explicitly it says the following.

\medskip

\noindent \textbf{Theorem \ref{main}} \emph{Let $V$ be a finitely generated pointwise finite dimensional $VIC(\mathbb{Z})$-module over $\mathbb{C}$,  $V$ admits a finite filtration $$V \supseteq V^1 \supseteq \dots \supseteq V^k = 0$$ by $VIC(\mathbb{Z})$-submodules such that $V / V^1$ is a torsion $VIC(\mathbb{Z})$-module (meaning $(V / V^1)_n = 0$ for all $n \gg 0$) and for each other successive subquotient $V^i / V^{i+1}$ there exists a bipartition $(\lambda^+_i, \lambda^-_i)$ such that $$(V^i / V^{i+1})_n = V_n(\lambda^+_i, \lambda^- _i) \otimes M^i_n$$ 
where the collection of $M^i_n$'s form a finitely generated $VIC(\mathbb{Z} / \ell_i \mathbb{Z})$-module for some integer $\ell_i$.}

\smallskip
\noindent \textbf{Remark:} Here $V_n(\lambda^+_i, \lambda^- _i)$ will denote a certain irreducible representation of $GL_n(\mathbb{C})$, which will naturally act on these spaces.

\vspace{.5cm}

\noindent The paper is structured as follows:

\begin{itemize}

\item Section \ref{fixedn} is primarily expository and discusses the finite dimensional representations of $SL_n(\mathbb{Z})$ for fixed $n \ge 3$ and the characterization of the irreducible representations via superrigidity and the congruence subgroup property.  Using this we introduce a notion of depth of an $SL_n(\mathbb{Z})$ representation, and use this to form effective improvements on the statement of superrigidity.

\item Section \ref{vicsect} contains the main results of the paper about finitely generated pointwise finite dimensional $VIC(\mathbb{Z})$-modules.  Other than the Noetherianity statement mentioned earlier and the main theorem stated above, this section also contains a technical result about extending from $VIC^\mathfrak{U}(R)$ to $VIC(R)$  which may be of interest to those in the field,.

\item Section \ref{applications} is about applications and extensions of the work here.  It contains a reformulation of the main result as an infinite rank superrigidity theorem, some results about $VIC(\mathbb{Z})$-modules over other rings including a Noetherianity statement in positive characteristic, and a characterization of $VIC(\mathbb{Z})$-modules of slow dimension growth.

\item Section \ref{future} is about open questions and directions of research coming out of this work. 

\end{itemize}

\section*{Acknowledgments}

Special thanks to Benson Farb for explaining superrigidity to me and to Andy Putman for initially suggesting looking at the depth of $VIC(\mathbb{Z})$-modules, as well as for a number of helpful discussions with both of them along the way.  Thanks to Peter Patzt, Steven Sam, and Nir Gadish for helpful conversations about $VIC$-modules while I was working on this. Thanks to Peter Patzt, Andrew Snowden, and Benson Farb for many helpful comments on an earlier draft of this paper. I'd also like to in general thank the UChicago geometry, topology, and representation theory working group, this work largely grew out of discussions held in the working group lounge.  This work was partially supported by NSF postdoctoral fellowship award 1703942.

\section{Effective superrigidity and depth for $SL_n(\mathbb{Z})$}\label{fixedn}

Before we start talking about compatible sequences of representations and $VIC(\mathbb{Z})$-modules we need to develop a bit of theory and language to talk about representations of a single $SL_n(\mathbb{Z})$. This part of the paper is mostly expository about the finite dimensional complex representation theory of $SL_n(\mathbb{Z})$. This section contains no new technical results, but we do state a number of consequences of existing results that appear to be new and possibly of independent interest from the rest of the paper.

 In particular, we use a characterization of finite dimensional representations of $SL_n(\mathbb{Z})$ via superrigidity and the congruence subgroup property to define a notion of \emph{depth} of a representation, which in a sense measures the failure of the representation to extend to $SL_n(\mathbb{C})$.  We give internal characterizations of depth in terms of the actions of certain matrices and subgroups, and bounds relating the dimension of a representation and its depth.

We'd like to think of this as a sort of effective version of superrigidity: Superrigidity for $SL_n(\mathbb{Z})$ tells us that any representation of $SL_n(\mathbb{Z})$ agrees with a representation of $SL_n(\mathbb{C})$ on \emph{some} finite index subgroup.  Here we are looking at properties of the representation to extract arithmetic information about \emph{which} finite index subgroup arises this way.

\subsection{Classification of irreducibles via superrigidity}

First let's review the classification of irreducible finite dimensional representations of $SL_n(\mathbb{Z})$ via superrigidity.

\subsubsection*{Algebraic representations}
One way to construct representations of  $SL_n(\mathbb{Z})$ is to take an algebraic (a.k.a. rational) representation $\phi : SL_n(\mathbb{C}) \rightarrow GL(V)$ and then obtain a restricted representation for $SL_n(\mathbb{Z})$ via the composition $$SL_n(\mathbb{Z}) \hookrightarrow SL_n(\mathbb{C}) \xrightarrow[]{\phi} GL(V)$$  We call such representations \emph{algebraic}.  The Borel density theorem (see \cite{Bor}) says that $SL_n(\mathbb{Z})$ is Zariski dense in $SL_n(\mathbb{C})$, this implies that if we restrict an irreducible algebraic representation to $SL_n(\mathbb{Z})$ it remains irreducible, and moreover a homomorphism of algebraic representations as $SL_n(\mathbb{Z})$-representations is also a homomorphism of $SL_n(\mathbb{C})$ representations.

\subsubsection*{Finite type representations}

Another way to construct representations of $SL_n(\mathbb{Z})$ is to take a representation $\psi: SL_n(\mathbb{Z} / \ell \mathbb{Z}) \rightarrow GL(W)$ and obtain a representation of $SL_n(\mathbb{Z})$ via the composition $$SL_n(\mathbb{Z}) \twoheadrightarrow SL_n(\mathbb{Z} / \ell \mathbb{Z}) \xrightarrow[]{\psi} GL(W)$$
We say that such representations are \emph{finite type}.  The congruence subgroup property for $SL_n(\mathbb{Z})$ ($n>2$), due independently to Mennicke and Bass-Lazard-Serre (\cite{Me} and \cite{BLS}), says that any homomorphism from $SL_n(\mathbb{Z})$ to a finite group factors through $SL_n(\mathbb{Z} / \ell \mathbb{Z})$ for some $\ell$. 

 Rephrasing this, a representation of finite type is one which factors through a finite quotient, which is the same thing as being a continuous representation of the profinite completion. The congruence subgroup property then tells us that the profinite completion of $SL_n(\mathbb{Z})$ is $\displaystyle \prod_p SL_n(\mathbb{Z}_p)$.  

\subsubsection*{Superrigidity and the full classification of irreducibles}

We can then mix these two types together and obtain more complicated representations by taking a tensor product $V \otimes W$ where $V$ is algebraic and $W$ is of finite type, and then taking direct sums of representations of this form.  Superrigidity for $SL_n(\mathbb{Z})$ (for $n \ge 3$) tells this is all we can do.  More precisely:
\begin{enumerate}

\item If $V$ is algebraic and $W$ is finite type and both are irreducible, then $V \otimes W$ is irreducible.

\item Two irreducible representations $V \otimes W$ and $V' \otimes W'$ of this form are isomorphic if and only if $V \cong V'$ and $W \cong W'$.

\item Every irreducible finite dimensional representation of $SL_n(\mathbb{Z})$ is of this form.

\item Every finite dimensional representation decomposes into a direct sum of irreducibles. 

\end{enumerate}

Note that if we take $V$ or $W$ to be trivial this in description we recover the finite type and algebraic cases, so they are included in this characterization.  For such a tensor product $M = V \otimes W$ as above we say that $V$ is the \emph{algebraic part} of $M$ and $W$ is the \emph{finite type part}.

\medskip

\noindent \textbf{Remark:} Superrigidity for $SL_n(\mathbb{Z})$ was proved by Bass, Milnor, and Serre (\cite{BMS} Theorem 16.2), however this explicit characterization of the irreducible representations given above seems hard to find spelled out in the literature.  In any case it follows easily from their formulation, and is certainly implicit in the computation of the proalgebraic completion of $SL_n(\mathbb{Z})$ by Bass et al. (in \cite{BLMM}) described in the next section.

\subsection{Categorical formulation of superrigidity}\label{catsup}

Categorically superrigidity can be expressed as saying that for $n>2$ $$\Rep(SL_n(\mathbb{Z})) \cong \Rep^{alg}(SL_n(\mathbb{C})) \boxtimes \mathlarger{\boxtimes_p} \Rep^{sm} (SL_n(\mathbb{Z}_p))$$
where $\Rep(SL_n(\mathbb{Z}))$ denotes the category of finite dimensional representations of $SL_n(\mathbb{Z})$, $\Rep^{alg}(SL_n(\mathbb{C}))$ denotes the category of finite dimensional algebraic representations of $SL_n(\mathbb{C})$, $\Rep^{sm} (SL_n(\mathbb{Z}_p))$ denotes the category of smooth finite dimensional representations of $SL_n(\mathbb{Z}_p)$ (i.e. those factoring through $SL_n(\mathbb{Z}/p^k)$ for some $k$), and $\boxtimes$ denotes the Deligne tensor product of tensor categories. 

  Explicitly since everything is semisimple this has objects given by direct sums of formal infinite products $V_{alg} \boxtimes V_2 \boxtimes V_3 \boxtimes V_5 \boxtimes V_7 \dots$ where all but finitely many of the terms $V_p \in \Rep^{sm} (SL_n(\mathbb{Z}_p))$ are the unit object (i.e. the trivial representation).  Categories obtained as such a tensor product satisfy a universal property categorifying the universal property of a tensor product of rings.
  
  The right hand side of the above equivalence should be thought of as the category of algebraic representations of the proalgebraic group $$SL_n(\mathbb{C}) \times \prod_p SL_n(\mathbb{Z}_p)$$ which is called the proalgebraic completion of $SL_n(\mathbb{Z})$ for this reason. Abstractly this means we could have recovered the proalgebraic group $SL_n(\mathbb{C}) \times \prod_p SL_n(\mathbb{Z}_p)$ by applying Tannakian reconstruction to the category of finite dimensional representations of $SL_n(\mathbb{Z})$.  More explicitly though, this means that the functor $$\Rep^{alg}(SL_n(\mathbb{C}) \times \prod_p SL_n(\mathbb{Z}_p)) \to \Rep(SL_n(\mathbb{Z}))$$ given by restriction to the ``diagonal" copy of $SL_n(\mathbb{Z})$ is an equivalence of categories.
  
These equivalences commute with all of the basic operations one might want to do with representations such as taking duals or tensor products. For example, the duality endofunctor on $\Rep(SL_n(\mathbb{Z}))$ clearly gets sent to the tensor product of duality endofunctors on the right hand side since if $M \cong V \otimes W$ where $V$ is algebraic and $W$ is finite type then $M^* \cong V^* \otimes W^*$. More categorically, this means the above restriction is a symmetric monoidal equivalence of Tannakian categories.

For our purposes it will be particularly important to note that these equivalences are also compatible with restriction of representations to Levi factors.  This is clear though from the above description of the equivalence coming from restriction:  If we want to restrict from $SL_{n+m}(\mathbb{C}) \times \prod_p SL_{n+m}(\mathbb{Z}_p)$ to $SL_n(\mathbb{Z}) \times SL_m(\mathbb{Z})$ then it doesn't matter if as an intermediate step we first restrict to $SL_{n+m}(\mathbb{Z})$ or to $SL_{n}(\mathbb{C}) \times \prod_p SL_{n}(\mathbb{Z}_p) \times SL_{m}(\mathbb{C}) \times \prod_p SL_{m}(\mathbb{Z}_p)$.

The theory of proalgebraic completions of superrigid groups was worked out by Bass et al. in \cite{BLMM} with the example of $SL_n(\mathbb{Z})$ appearing explicitly in Section 6.8. For general background on the Deligne tensor product of categories, tensor categories, and Tannakian reconstruction we refer to \cite{EGNO}.

\subsection{Depth of representations}

If $V$ is a representation of $SL_n(\mathbb{Z})$ we define the \emph{depth} of $V$ to be the smallest $\ell$ such that the representation is algebraic when restricted to the congruence subgroup $\Gamma_n(\ell)$.  The $p$-\emph{depth} of $V$ is the largest power $k$ such that $p^k$ divides the depth. 

If $V$ is irreducible then the depth is just the smallest $\ell$ such that we can write $V$ as a tensor product $V_{alg} \otimes V_{fin}$ with $V_{alg}$ algebraic and $V_{fin}$ is a finite type representation factoring through $SL_n(\mathbb{Z}/\ell\mathbb{Z})$. For $V$ reducible the depth is just the least common multiple of the depths of its irreducible components.

\

\noindent A few easy observations about this:

\begin{enumerate}

\item $V$ has depth $1$ if and only if $V$ is algebraic.

\item $V^*$ has the same depth as $V$.

\item The depth of $V \otimes W$ and $V\oplus W$ are at most the least common multiple of the depths of $V$ and $W$. 

\item If $n>3$ and we restrict $V$ from $SL_n(\mathbb{Z})$ to $SL_{n-1}(\mathbb{Z})$ the depth remains the same.

\end{enumerate}

The category of representations of depth dividing $\ell$ is a Tannakian category equivalent to the category of algebraic representations of the (genuine, not pro-) algebraic group $SL_n(\mathbb{C}) \times SL_n(\mathbb{Z}/\ell \mathbb{Z})$.  The full category of representations is just a direct limit (i.e. a union) of these subcategories.

\bigskip

We will be interested in how to understand and bound the depth of representations, and especially families of representations.  If we have a sequence of representations $V_3,V_4,V_5, \dots$ where each $V_k$ is a representation of $SL_k(\mathbb{Z})$ (say coming from a VIC-module), we say that the sequence has $\emph{bounded depth}$ if there is some $\ell$ such that each $V_i$ has depth dividing $\ell$.  In this case we will refer to the smallest such $\ell$ as the depth of the sequence of representations.  We will also say that such a sequence is \emph{eventually algebraic} if there is some $N$ such that $V_n$ is algebraic for all $n > N$.

\subsection{Internal characterizations of depth}

Right now, the depth of $V$ is defined in terms of how we can write $V$ as a virtually algebraic representation.  While we know this is always possible by superrigidity, it's not clear how we can actually figure out the depth of any particular representation $V$. We'll now give a different characterization of the depth that is more intrinsic to the representation itself.

\subsubsection*{Depth and the action of elementary matrices}

Let $E =  E_{1,2} \in SL_n(\mathbb{Z})$ be the elementary matrix where $a_{i,j} = 1$ if $i = j$ or if $(i,j) = (1,2)$ and $0$ otherwise.  This is conjugate to any other elementary matrix, where the off-diagonal $1$ is somewhere else, so everything we say will also hold for any other elementary matrix. First we'd like to state an important fact about $E$ as it relates to the congruence subgroups $\Gamma_n(\ell)$.  

\begin{lemma}\label{normgen}
 $\Gamma_n(\ell)$ is the smallest normal subgroup of $SL_n(\mathbb{Z})$ containing $E^\ell$. 
\end{lemma}

This explicit statement was formulated by Mennicke, and is the main theorem of \cite{Me}. Bass-Lazard-Serre also proved an equivalent statement in \cite{BLS}, and it is the main technical ingredient in the proof of the congurence subgroup property for $SL_n(\mathbb{Z})$.  With this, we are ready to give the following characterization of depth in terms of how $E$ acts:

\begin{proposition} If $V$ is a finite dimensional representation of $SL_n(\mathbb{Z})$ then every eigenvalue of $E$ is a root of unity, and 

\begin{enumerate}

\item The depth of $V$ is the least common multiple of the orders of these eigenvalues.

\item If $V$ is irreducible then the depth of $V$ is the largest order of these eigenvalues.

\end{enumerate}

\end{proposition}
Note: In general $E$ does not act by a semisimple matrix. In fact, $E$ acts by a semisimple matrix if and only if $V$ is a finite type representation. In particular, a representation of $SL_n(\mathbb{Z})$ is unitary if and only if it is finite type.

\medskip

\noindent \textbf{Proof:}  First note that $E$ acts unipotently in any algebraic representation, so in general the only eigenvalues that occur are those that appear in finite type representations. Then by the Chinese remainder theorem it is enough to check part 2 in the case where $\ell = p^k$ is a prime power, and we are viewing $E$ as an element of $SL_n(\mathbb{Z}/p^k\mathbb{Z})$.

 Here $E$ has order $p^k$, so all of its eigenvalues have orders dividing that automatically. Moreover if all of the eigenvalues were have strictly smaller order then $E^{p^{k-1}}$ acts trivially so therefore by Lemma \ref{normgen} the action factors through $SL_n(\mathbb{Z}/p^{k-1}\mathbb{Z})$ and thus has smaller depth. $\square$
 
 \medskip
 
Instead, we could state this just in terms of how an $\ell$th power of an elementary matrix acts. Again these matrices are all conjugate so we'll focus on $E^\ell =  E_{1,2}^\ell$ to state the following corollary:

\begin{corollary} \label{elemdepth} A representation of $SL_n(\mathbb{Z})$ has depth dividing $\ell$ if and only if $E^\ell$ acts unipotently.

\end{corollary}

\subsubsection*{Action of the unipotent subgroup}

We just looked at how a single elementary matrix acted on representations, now let's extend that analysis to the entire group of unipotent uppertriangular matrices. Let $U_n(R) \subset SL_n(R)$ denote the subgroup of upper triangular matrices with ones on the diagonal, for any ring $R$.  Let's recall some facts about these groups and how they act on familiar (i.e. algebraic and finite type) representations:

\begin{itemize}
\item Algebraic representations of $U_n(\mathbb{C})$ do not in general decompose as a direct sum of irreducibles, and in fact any algebraic representation of $U_n(\mathbb{C})$ is an iterated extension of trivial representations. 

\item $U_n(\mathbb{Z})$ is Zariski dense in $U_n(\mathbb{C})$, so any $U_n(\mathbb{Z})$-invariant subspace of an algebraic representation is also $U_n(\mathbb{C})$-invariant.

\item In an algebraic representation, any vector which is not $U_n(\mathbb{C})$ invariant has an infinite (and unbounded in the usual topology) $U_n(\mathbb{Z})$-orbit.

\item Any irreducible finite dimensional representation of $SL_n(\mathbb{C})$ has a unique $U_n(\mathbb{C})$ invariant line spanned by a highest weight vector for the diagonal torus.  

\item Of course $U_n(\mathbb{Z}/\ell \mathbb{Z})$ is a finite group and therefore any representation of it decomposes as a direct sum of irreducible representations. Also the orbit of any vector is obviously finite.

\item If the prime factorization of $\ell$ is $p_1^{k_1} p_2^{k_2} \dots p_m^{k_m}$ then $$U_n(\mathbb{Z}/\ell \mathbb{Z}) \cong U_n(\mathbb{Z}/p_1^{k_1} \mathbb{Z}) \times  U_n(\mathbb{Z}/p_2^{k_2} \mathbb{Z}) \times \dots \times U_n(\mathbb{Z}/p_m^{k_m} \mathbb{Z})$$ and each $U_n(\mathbb{Z}/p_i^{k_i} \mathbb{Z})$ is a $p_i$-group. In particular $p_1, p_2, \dots$ and $p_m$ are the only prime numbers dividing the dimensions of irreducible representations of $U_n(\mathbb{Z}/\ell \mathbb{Z})$. 

\end{itemize}

Note that we don't expect any version of the congruence subgroup property to hold for $U_n(\mathbb{Z})$, so the full representation theory of $U_n(\mathbb{Z})$ may be difficult to describe. However those representations that appear in restrictions of finite dimensional $SL_n(\mathbb{Z})$ are fairly restricted and inherit some nice properties for free. Combining some of these facts above gives us a bit of structure theory for the action of $U_n(\mathbb{Z})$ on $SL_n(\mathbb{Z})$-representations.

\begin{corollary} Let $V = V_{alg}\otimes V_{fin}$ be a finite dimensional irreducible representation of $SL_n(\mathbb{Z})$ with $V_{alg}$ and $V_{fin}$ its algebraic and finite type parts.  Then as a $U_n(\mathbb{Z})$ representation:

\begin{itemize}

\item The socle of $V$ is isomorphic to $V_{fin}$ as $U_n(\mathbb{Z})$-representations. In particular, this action on the socle factors through $U_n(\mathbb{Z}/\ell \mathbb{Z})$ where $\ell$ is the depth of $V$. 

\item This socle can alternatively be characterized as the space of all vectors $\vec{v} \in V$ such that the $U_n(\mathbb{Z})$-orbit of $\vec{v}$ is finite.

\item  $V$ has a filtration by $U_n(\mathbb{Z})$-invariant subspaces $0 =V_0 \subset V_1 \subset ... \subset V_k = V$ such that each quotient $V_i / V_{i-1}$ is isomorphic to this socle and $k$ is equal to the dimension of $V_{alg}$.

\end{itemize}

\end{corollary}

Next define $\Gamma U_n (\ell)$, the upper triangular congruence subgroup of level $\ell$, to be the kernel of the map from $U_n(\mathbb{Z})$ to $U_n(\mathbb{Z}/\ell \mathbb{Z})$.  Explicitly, this is just the group of upper triangular $n\times n$ integer matrices with ones on the diagonal and $\ell$ dividing every entry off the diagonal. Here are a few easy observations about the action of $\Gamma U_n (\ell)$ on representations of $SL_n(\mathbb{Z})$.

\begin{itemize}

\item $\Gamma U_n (\ell)$ is still Zariski dense in $U_n(\mathbb{C})$ so any $\Gamma U_n (\ell)$-invariant subspace of an algebraic representation is also $U_n(\mathbb{C})$ invariant. In particular an irreducible algebraic representation of $SL_n(\mathbb{C})$ has a unique $\Gamma U_n (\ell)$ invariant vector.

\item $\Gamma U_n (\ell)$ acts trivially on any representation of $SL_n(\mathbb{Z})$ which factors through $SL_n(\mathbb{Z}/\ell \mathbb{Z})$.

\item  $E^\ell \in \Gamma U_n (\ell)$ and if therefore by Corollary \ref{elemdepth} if $V$ has depth not dividing $\ell$ then there is an element of $\Gamma U_n (\ell)$ acting non-unipotently.

\end{itemize}
In particular this gives the following extension of Corollary \ref{elemdepth}

\begin{corollary}
A representation of $SL_n(\mathbb{Z})$ has level dividing $\ell$ if and only if every element of $\Gamma U_n (\ell)$ acts unipotently.

\end{corollary}

\subsection{Bounding depth in terms of dimension}

The following observation is what started this entire investigation:  $SL_n(\mathbb{Z})$ has an obvious nontrivial algebraic representation of dimension $n$, but the non-trivial smallest finite type representation the author knew how to write down (functions on $\mathbb{P(}\mathbb{F}_2^n)$ with total sum $0$) has dimension $2^n - 2$. 

 This suggested that perhaps all low dimensional representation had to be algebraic, something that could have some interesting applications to representation stability.  Indeed this turns out to be the case, and more generally we can bound the depth of a representation just in terms of the dimension.

Ultimately this is just a statement about the minimal dimension of a representation of $SL_n(\mathbb{Z}/\ell \mathbb{Z})$ which does not factor through some $SL_n(\mathbb{Z}/\ell' \mathbb{Z})$ with $\ell ' < \ell$.  We'll now summarize what is known about these minimal dimensions. First up we'll mention the case where $\ell = p$ is a prime number, this is the most well studied case and here we have a precise answer due to Tiep and Zalesski:

\begin{theorem}{\textbf{(\cite{TZ} Theorem 1.1)}} In the case where $\ell = p$ is prime and if $n \ge 5$ or $p > 3$, the smallest dimension of a nontrivial representation of $SL_n(\mathbb{Z}/p \mathbb{Z})$ is $\frac{p^n - p}{p-1} = p^{n-1} + p^{n-2} + \dots + p$.  Moreover the only representation having this dimension is the space of complex valued functions on the finite set $\mathbb{P(}\mathbb{F}_p^n)$ with total sum $0$.
\end{theorem}

Next up is the case where $\ell = p^k$ is a prime power.  Here we don't have exact values but we some bounds due to Bardestani and Mallahi-Karai:

\begin{theorem}\label{pkbound}{\textbf{(\cite{BMK} Theorem 1)}} Any representation of $SL_n(\mathbb{Z}/p^k \mathbb{Z})$ which does not factor through $SL_n(\mathbb{Z}/p^{k-1} \mathbb{Z})$ has dimension at least $$(p^k - p^{k-1})p^{(n-2)k} = p^{(n-1)k} (1- \frac{1}{p})$$
\end{theorem}

It's worth noting that while this isn't an exact answer, it is approximately the size of the finite projective space $(\mathbb{Z}/p^k \mathbb{Z})^n / (\mathbb{Z}/p^k \mathbb{Z})^\times$, so this bound can't be significantly improved. For our purposes we are mostly interested in the rate of growth of these dimensions as a function of $n, k,$ and $p$, so these approximate answers are enough. 

However it is reasonable to conjecture that the space of complex functions with total sum $0$ on this finite projective space gives the unique smallest such representation (with perhaps finitely many exceptions), but we don't pursue this direction any further.

\medskip

\noindent \textbf{Remark:} We'll note that Bardestani and Mallahi-Karai actually state this result as a bound on the dimension of a \emph{faithful} representation of $SL_n(\mathbb{Z}/p^k \mathbb{Z})$, however their proof actually proves this slightly stronger statement.  Explicitly, faithfulness only essentially appears in their proof of ``Lemma 3", where they arrive at a contradiction of the faithfulness of a representation by showing that an elementary matrix acts with order dividing $p^{k-1}$. However in fact this implies not just that the action isn't faithful, but that it factors through $SL_n(\mathbb{Z}/p^{k-1} \mathbb{Z})$ by Lemma \ref{normgen}.

\medskip

For general $\ell$ the Chinese remainder theorem tells us that if the prime factorization of $\ell$ is $p_1^{k_1} p_2^{k_2} \dots p_m^{k_m}$ then $$SL_n(\mathbb{Z}/ \ell \mathbb{Z}) = SL_n(\mathbb{Z}/p_1^{k_1} \mathbb{Z}) \times SL_n(\mathbb{Z}/p_2^{k_2} \mathbb{Z}) \times \dots \times SL_n(\mathbb{Z}/p_m^{k_m} \mathbb{Z})$$
Therefore the smallest irreducible representation of $SL_n(\mathbb{Z}/ \ell \mathbb{Z})$ not factoring through some $SL_n(\mathbb{Z}/ \ell' \mathbb{Z})$ with $\ell' < \ell$ is just a tensor product of the smallest irreducible representations of each $SL_n(\mathbb{Z}/p_i^{k_i}\mathbb{Z})$ not factoring through $SL_n(\mathbb{Z}/p_i^{k_i - 1}\mathbb{Z})$.  We obtain the following corollary:

\begin{corollary}\label{dimdepth} An irreducible representation of $SL_n(\mathbb{Z})$ of depth $\ell$ has dimension at least $$\ell^{n-1} \cdot \prod_{p | \ell} (1 - \frac{1}{p}) \ge \ell^{n-2}$$
\end{corollary} 
(Writing it as $\ell^{n-2}$ is just for convenience, since we will be mostly interested in the rate of growth in $n$ anyway. )

\medskip

\noindent \textbf{Remark:} We'll note that if we drop the irreducible condition this corollary is false as stated. As defined, the smallest depth $\ell$ representation is actually the direct sum of the smallest irreducible representations of each $SL_n(\mathbb{Z}/p_i^{k_i}\mathbb{Z})$ not factoring through $SL_n(\mathbb{Z}/p_i^{k_i - 1}\mathbb{Z})$ rather than the tensor product.

\medskip

 A slightly improved special case of this, obtained using the exact answers from Tiep and Zalesskii is the following corollary:

\begin{corollary}\label{lowdimalg}  For $n \ge 5$ any representation of $SL_n(\mathbb{Z})$ of dimension less than $2^n - 2$ extends to an algebraic representation of $SL_n(\mathbb{C})$.

\end{corollary}

One may think of Corollaries \ref{dimdepth} and \ref{lowdimalg} as effective strengthenings of superrigidity.  We knew from superrigidity that any homomorphism $SL_n(\mathbb{Z}) \to GL_N(\mathbb{C})$ agrees with an algebraic representation of $SL_n(\mathbb{C})$ along some congruence subgroup $\Gamma_n(\ell)$. These corollaries tell us that  in fact we can bound this $\ell$ in terms of $N$, the dimension of the representation.

\medskip

The bulk of the remainder of the paper will be about $VIC(\mathbb{Z})$-modules, which are certain sequences of representations with compatibility constraints between them.  First though we'd like to state a quick corollary about sequences of finite dimensional representations $V_3,V_4,V_5, \dots$ where each $V_k$ is a representation of $SL_k(\mathbb{Z})$ but without any compatibility assumptions between the representations.

 Say the sequence has \emph{exponential growth} if there is some $C$ such that $\text{dim}(V_n) \le C^n$ for all $n$.  Similarly we say the sequence has \emph{polynomial growth} if there is some $k$ such that $\text{dim}(V_n) \le n^k$.  We have the following corollaries:

\begin{corollary} Let $V_* = V_3,V_4,V_5, \dots$  be a sequence of representations where each $V_k$ is a representation of $SL_k(\mathbb{Z})$.

\begin{enumerate}

\item If $V_*$ has exponential growth then it has bounded depth.

\item If $V_*$ has polynomial growth then it is eventually algebraic.

\end{enumerate}

\end{corollary}

\section{Pointwise finite dimensional $VIC(\mathbb{Z})$-modules}\label{vicsect}

Suppose $R$ be a commutative ring. If $M$ and $N$ are free $R$-modules, a \emph{splitable injection} from $M$ to $N$ is a linear injection $f: M \to N$ which admits a linear left inverse $g: N \to M$  satisfying $g\circ f (m) = m \ \forall m \in M$. A \emph{split injection} from $M$ to $N$ is the data of a splitable linear injection $f: M \to N$ along with a choice of such a left inverse.

Note that the data of a linear left inverse to an injection $f: M \to N$ is equivalent to the data of a complementary subspace $M^\perp$ so that $N = \text{im}(f) \oplus M^\perp$.  The equivalence being that we can take the complementary subspace to be the kernel of the left inverse map.

$VIC(R)$ is the category of free $R$-modules of rank at least 3 (see warning below) with split injections between them. A $VIC(R)$-module over a fixed ground ring $k$ is a functor from $VIC(R)$ to the category of modules over $k$. We will restrict to the case when $k = \mathbb{C}$.

A $VIC(R)$-module $F:VIC(R) \to \text{Mod}_k$ is \emph{finitely generated} if there are finitely many vectors $v_1, v_2, \dots, v_k$ in $F(X_1), F(X_2), \dots F(X_k)$ such that any vector in any $F(X)$ can be written as a linear combination of images of the $v_i$'s under various split inclusions from $X_i$ into $X$.

\medskip

\noindent \textbf{Warning:} The rank at least 3 condition in our definition is non-standard.  It is harmless for the type of questions we will care about though. Any $VIC(R)$-module in the standard sense restricts to one in this sense, and every $VIC(R)$-module in this sense extends to one in the usual sense (although not uniquely). Mostly for us it just means we don't have to modify our statements to worry about what happens for $GL_2(\mathbb{Z})$ where superrigidity and the congruence subgroup property fail.

\subsubsection{An alternative definition and $VIC^{\mathfrak{U}}(R)$-modules}

While this language of functors from a category into vector spaces is extremely useful for applications, it may not be the clearest for understanding what $VIC(R)$ modules actually look like from a representation theoretic perspective. 

 Instead it is often easier to work with an alternative characterization of $VIC(R)$-modules due to Randal-Williams and Wahl (\cite{RWW} Prop. 4.2).   A $VIC(R)$-module $V$ is equivalent to the data of a sequence of vector spaces $V_3, V_4, V_5, \dots$ over $k$ such that:

\begin{enumerate}

\item Each $V_n$ is a representation of $GL_n(R)$.

\item For each $n$ there is a linear map $V_n \to V_{n+1}$ which is $GL_n(R)$-equivariant for the standard inclusion (i.e. in the top left corner) of $GL_n(R)$ into $GL_{n+1}(R)$.
 
\item The standard complementary copy of $GL_m(R)$ to $GL_n(R)$ inside $GL_{n+m}(R)$ (i.e. in the bottom right corner) acts trivially on the image of $V_n$ inside $V_{n+m}$.
 
\end{enumerate}

This reformulation will be convenient for defining an important generalization of $VIC$-modules.  Let $\mathfrak{U} 
\subseteq R^\times$ be a subgroup of units in $R$, and let $SL^\mathfrak{U}_n(R)$ denote the group of $n \times n$  matrices with coefficients in $R$ and determinant in $\mathfrak{U}$.  A $VIC^{\mathfrak{U}}(R)$-module is the data of a sequence of vector spaces $V_3, V_4, V_5, \dots$ over $k$ such that:

\begin{enumerate}

\item Each $V_n$ is a representation of $SL^{\mathfrak{U}}_n(R)$.

\item For each $n$ there is a linear map $V_n \to V_{n+1}$ which is $SL^{\mathfrak{U}}_n(R)$-equivariant equivariant for the standard inclusion of $SL^{\mathfrak{U}}_n(R)$ into $SL^{\mathfrak{U}}_{n+1}(R)$.
 
\item The complementary copy of $SL^{\mathfrak{U}}_m(R)$ acts trivially on the image of $V_n$ inside $V_{n+m}$.
 
\end{enumerate}

Most often we will care about the case when $\mathfrak{U} = \{ -1, 1 \}$ in this case we will use the notation $SL^\pm(R)$ and $VIC^\pm(R)$.  In particular this will arise since $1$ and $-1$ are the only units in $\mathbb{Z}$, which implies that $SL_n^\pm(\mathbb{C})$ is the Zariski closure of $GL_n(\mathbb{Z})$ in $GL_n(\mathbb{C})$, and that $SL_n^\pm(\mathbb{Z}/\ell \mathbb{Z})$ is the image of $GL_n(\mathbb{Z})$ in $GL_n(\mathbb{Z}/\ell \mathbb{Z})$.  Randal-Williams and Wahl's alternative characterization of $VIC$-modules, along with our definition of $VIC^{\mathfrak{U}}(R)$-modules gives us the following corollary:

\begin{corollary}\label{pointwise} \
\begin{enumerate}

\item If a $VIC(\mathbb{Z})$-module $V$ has the property that for some fixed $\ell \in \mathbb{N}$ the action of $GL_n(\mathbb{Z})$ on $V_n$ factors through $VIC(\mathbb{Z}/\ell \mathbb{Z})$ for each $n$, then it has the structure of a $VIC^\pm(\mathbb{Z}/\ell \mathbb{Z})$-module. Moreover, it is finitely generated as a $VIC^\pm(\mathbb{Z}/\ell \mathbb{Z})$-module if and only if it is finitely generated as a $VIC(\mathbb{Z})$-module.

\item  If a $VIC(\mathbb{Z})$-module $V$ has the property that each $V_n$ is an algebraic representation of $GL_n(\mathbb{Z})$, then it is the restriction a  $VIC^\pm(\mathbb{C})$-module. Moreover it is finitely generated as a $VIC^\pm(\mathbb{C})$-module if and only if it is finitely generated as a $VIC(\mathbb{Z})$-module.

\end{enumerate}
\end{corollary}

Often we will care about $VIC$-modules with additional constraints on the representations $V_n$.  We will use the term ``pointwise" to describe this situation.  For example a $VIC(\mathbb{Z})$-module is pointwise finite type if each $V_n$ is a representation of $GL_n(\mathbb{Z})$ factoring through a finite group.  Most importantly we will be interested in the case of pointwise finite dimensional $VIC(\mathbb{Z})$-modules, where each $V_n$ is a finite dimensional representation of $GL_n(\mathbb{Z})$.

\subsubsection*{Weak superrigidity for $VIC(\mathbb{Z})$}

The categorical formulation of superrigidity from Section \ref{catsup}, saying that the proalgebraic completion of $SL_n(\mathbb{Z})$ is $SL_n(\mathbb{C}) \times \prod_p SL_n(\mathbb{Z}_p)$, tells us that any finite dimensional $SL_n(\mathbb{Z})$-representation can be equipped with commuting actions of $SL_n(\mathbb{C})$ and $SL_n(\mathbb{Z}_p)$ for all $p$, such that the original action of $SL_n(\mathbb{Z})$ is given by the diagonal inclusion into $SL_n(\mathbb{C}) \times \prod_p SL_n(\mathbb{Z}_p)$. Combining this with formulation of $VIC$-modules due to Randal-Williams and Wahl a priori gives us the following weak version of superrigidity for $VIC(\mathbb{Z})$-modules:

\begin{proposition}\label{weaksr}

If $V$ is a pointwise finite dimensional $VIC(\mathbb{Z})$-module then the commuting actions of $SL^\pm_n(\mathbb{C})$ and $SL^\pm_n(\mathbb{Z}_p)$ on $V_n$ define the structures of a $VIC^\pm(\mathbb{C})$-module and a $VIC^\pm(\mathbb{Z}_p)$-module on the same underlying vector spaces.

\end{proposition}

\noindent \textbf{Proof:} We already have a sequence of vector spaces $V_n$ with actions of $SL^\pm_n(\mathbb{C})$ and $SL^\pm_n(\mathbb{Z}_p)$ as well as inclusions $V_n \to V_{n+1}$ which are $GL_n(\mathbb{Z})$-equivariant. The functorial nature of superrigidity implies these linear maps are also $SL^\pm_n(\mathbb{C})$ and $SL^\pm_n(\mathbb{Z}_p)$ equivariant, so we just need to check that the image of $V_n$ in $V_{n+m}$ is fixed by the actions of $SL^\pm_m(\mathbb{C})$ and $SL^\pm_m(\mathbb{Z}_p)$. This is true since in general any vector that is $GL_m(\mathbb{Z})$-invariant is also invariant for the induced actions $SL^\pm_m(\mathbb{C})$ and $SL^\pm_m(\mathbb{Z}_p)$. $\square$

\medskip

\noindent \textbf{Warning:} If $V$ is a finitely generated pointwise finite dimensional $VIC(\mathbb{Z})$-module then in general $V$ will not be finitely generated for the induced $VIC^\pm(\mathbb{C})$-module and $VIC^\pm(\mathbb{Z}_p)$-module structures.  This should not be that surprising though, note that even a finitely generated $VIC(\mathbb{Z}/pq\mathbb{Z})$-module is usually not finitely generated as a $VIC(\mathbb{Z}/p\mathbb{Z})$ or $VIC(\mathbb{Z}/q\mathbb{Z})$-module.

\subsubsection*{Tensor products of $VIC$-modules}

Suppose $V_n$ and $W_n$ are $VIC(R)$-modules, we can consider the pointwise tensor product $(V \otimes W)_n = V_n \otimes W_n$.  This naturally inherits the structure of a $VIC(R)$-module.  Miller and Wilson (\cite{MW} Prop 4.6) have shown that if $R$ is a PID then if $V$ and $W$ are $VIC(R)$-modules generated in finite degree, then their tensor product is also generated in finite degree\footnote{They mostly consider the case where $R$ is a field, and obtain explicit (and tight) bounds on the degree of generation in this case. However, as they note just after their proof of Lemma 4.4, the same proof applies to the case when $R$ is a PID, but you get weaker degree bounds.}.

While in general being generated in finite degree is a weaker condition than being finitely generated, we'll note that in the case when all of the representations are finite dimensional these two notions coincide.  In particular in our context we have that:

\begin{corollary} \label{tensor}

If $V$ and $W$ are finitely generated, pointwise finite dimensional $VIC(\mathbb{Z})$-modules, then $V \otimes W$ is also finitely generated and pointwise finite dimensional.

\end{corollary}

\subsubsection*{Shifts of $VIC$-modules}

Given a $VIC(\mathbb{Z})$-module $V: VIC(\mathbb{Z}) \to \text{Vec}_\mathbb{C}$ we may perform the shift operation to obtain a new $VIC(\mathbb{Z})$-module $SV$ defined functorially by $SV(M) = V(M\oplus \mathbb{Z})$  and $SV(f,C) = V(f\oplus \text{Id}, C)$. In terms of the underlying $GL_n(\mathbb{Z})$-representations this looks like $$SV_n = \text{Res}^{GL_{n+1}(\mathbb{Z})}_{GL_n(\mathbb{Z})}(V_{n+1}).$$

It follows from Proposition 3.27 in \cite{MPW} that if $V$ is generated in degrees $\le d$ then $SV$ is generated in degrees $\le 2d+1$ (although we suspect that this $2d+1$ can be improved to $d$ provided $d\ge2$). In particular though, again since under the assumption that our $VIC(\mathbb{Z})$-modules are pointwise finite dimensional the notions of being finitely generated and of being generated in finite degree coincide, we have the following corollary:

\begin{corollary}\label{shiftfg}

If $V$ is a finitely generated pointwise finite dimensional $VIC(\mathbb{Z})$-module, then $SV$ is also finitely generated and pointwise finite dimensional.
\end{corollary}

\subsubsection*{Stability of covariants}

Given a $VIC(\mathbb{Z})$-module $V$ we obtain a (single) map of covariants
$$(V_n)_{GL_n(\mathbb{Z})} \to (V_{n+1})_{GL_{n+1}(\mathbb{Z})} $$
for each $n \in \mathbb{N}$.  The collection of these spaces of covariants itself forms $VIC(\mathbb{Z})$-module -- it is the quotient by the submodule generated by all vectors of the form $\vec{v} - g\vec{v}$.  In particular if $V$ is finitely generated then so is any quotient of it, so this is a finitely generated $VIC(\mathbb{Z})$-module where each $V_n$ is a direct sum of trivial representations. 

Therefore for any finitely generated $VIC(\mathbb{Z})$-module the spaces of coinvariants are all finite dimensional and moreover the maps $(V_n)_{GL_n(\mathbb{Z})} \to (V_{n+1})_{GL_{n+1}(\mathbb{Z})} $ are isomorphisms for all $n \gg 0$, in particular the dimensions stabilize. Note that since finite dimensional representations of $GL_n(\mathbb{Z})$ over $\mathbb{C}$ are semisimple the covariants and invariants in a pointwise finite dimensional $VIC(\mathbb{Z})$-module can be identified pointwise, but in general the invariants do not form a $VIC(\mathbb{Z})$-submodule.

These spaces of covariants together with the maps between them can be thought of as a graded $\mathbb{C}[T]$-module $$\Phi_0(V) = \bigoplus_{n \in \mathbb{N}} (V_n)_{GL_n(\mathbb{Z})}$$
where $T$ acts by the maps $(V_n)_{GL_{n}(\mathbb{Z})} \to (V_{n+1})_{GL_{n+1}(\mathbb{Z})}$. The above reasoning tells us that if $V$ is finitely generated then $\Phi_0(V)$ is finitely generated over $\mathbb{C}[T]$ as well.

\subsubsection*{More stability of covariants}

Generalizing what we just did, if we fix an $a \in \mathbb{N}$ we may consider the $GL_{n-a}(\mathbb{Z})$ covariants inside $V_n$, for the inclusion of $GL_{n-a}(\mathbb{Z})$ in $GL_{n}(\mathbb{Z})$ as the last $(n-a)$ coordinates.  Our structure maps from being a $VIC(\mathbb{Z})$-module again descend to a map of covariants to give us a map $$(V_n)_{GL_{n-a}(\mathbb{Z})} \to (V_{n+1})_{GL_{n-a+1}(\mathbb{Z})} $$
which this time intertwines the actions of $GL_a(\mathbb{Z})$ acting on both factors via the standard embeddings of $GL_a(\mathbb{Z})$ into $GL_{n}(\mathbb{Z})$ and $GL_{n+1}(\mathbb{Z})$. 

 As before it is convenient to take the direct sum of these spaces and view it as a graded $\mathbb{C}[T]$-module (with a $GL_a(\mathbb{Z})$ action respecting the grading) defining:

$$\Phi_a(V) = \bigoplus_{n > a} (V_n)_{GL_{n-a}(\mathbb{Z})}$$ 
with $T$ acting via our maps $(V_n)_{GL_{n-a}(\mathbb{Z})} \to (V_{n+1})_{GL_{n-a+1}(\mathbb{Z})} $. However unlike when $a = 0$ we'll note that these covariant spaces are not naturally the quotient of $V$ by a $VIC(\mathbb{Z})$-submodule, and in general for $V$ finitely generated $\Phi_a(V)$ may not be finitely generated over $\mathbb{C}[T]$.

However the savvy reader may note that $\Phi_a(V) \cong \Phi_0(S^aV)$ where $S^a$ is the shift operator iterated $a$ times, so finite generation of $\Phi_a(V)$ does follow if we know $S^aV$ is finitely generated.  Combining this observation with Corollary \ref{shiftfg} we obtain the following corollary:

\begin{corollary}\label{Phia}

If $V$ is a finitely generated, pointwise finite dimensional $VIC(\mathbb{Z})$-module then $\Phi_a(V)$ is finitely generated over $\mathbb{C}[T]$, and the maps $(V_n)_{GL_{n-a}(\mathbb{Z})} \to (V_{n+1})_{GL_{n-a+1}(\mathbb{Z})} $ are isomorphisms for all $n \gg 0$.

\end{corollary}


\subsection{Noetherianity for finite dimensional $VIC(\mathbb{Z})$-modules}

Noetherianity results are the cornerstone of the theory of representation stability.  They are what allow us to perform spectral sequence arguments to prove finite generation of many of the most interesting and important instances of representation stability.

Putman and Sam proved a Noetherianity result for finitely generated $VIC(R)$-modules when $R$ is a finite ring, but showed that if $R$ is infinite (and in particular if $R = \mathbb{Z})$ that finitely generated $VIC(R)$-modules may have submodules which are not finitely generated. However we'll note that their counterexamples are all infinite dimensional. 

On the other hand, Patzt showed that if we consider the category of finitely generated $VIC(\mathbb{C})$-modules which are pointwise algebraic representations, then Noetherianity holds and any submodule of such a $VIC(\mathbb{C})$-module is also finitely generated.

Superrigidity tells us that finite dimensional representations of $SL_n(\mathbb{Z})$ are built out of algebraic representations of $SL_n(\mathbb{C})$ and representations of $SL_n(\mathbb{Z} / \ell \mathbb{Z})$, both of which have Noetherianity properties.  So one might expect a Noetherianity result to hold for pointwise finite dimensional $VIC(\mathbb{Z})$-modules and indeed that is the case. We have the following theorem:

\begin{theorem} \label{Noetherianity}
 If $V$ is a finitely generated, pointwise finite dimensional $VIC(\mathbb{Z})$-module then any submodule $W \subseteq V$ is also finitely generated.

\end{theorem}

\noindent \textbf{Proof:}  We will follow Patzt's proof of Noetherianity for pointwise algebraic $VIC(\mathbb{C})$-modules (\cite{Patzt1} Theorem C), which itself is based off the proof of Noetherianity for $FI$-modules in characteristic zero due to Church, Ellenberg, and Farb (\cite{CEF} Theorem 1.3).

Suppose $V$ is generated in degree $\le a$, then $\Phi_a(V)$ is a finitely generated graded $\mathbb{C}[T]$ module by Corollary \ref{Phia}.  Since $\mathbb{C}[T]$ is Noetherian this implies that the submodule $\Phi_a(W)$ is also finitely generated.

If we choose homogeneous generators $x_1, x_2, \dots x_r$ for $\Phi_a(W)$, with $x_i \in \Phi_a(W)_{n_i}$, then we can find lifts $w_1, w_2, \dots, w_r$ with $w_i \in W_{a + n_i}$ projecting to $x_i$.  Let $\widetilde{W} \subseteq W$ be the submodule generated by these $w_1, w_2, \dots, w_r$. It is enough to show that $\widetilde{W}_n = W_n$ for $n \gg 0$.

By construction we have that $\Phi_a(\widetilde{W}) = \Phi_a(W)$ so explicitly this means that $(\widetilde{W}_n)_{GL_{n-a}(\mathbb{Z})} = (W_n)_{GL_{n-a}(\mathbb{Z})} $ for all $n$.  So $(W_n / \widetilde{W}_n)_{GL_{n-a}(\mathbb{Z})} = 0$ and by semisimplicity therefore $(W_n / \widetilde{W}_n)^{GL_{n-a}(\mathbb{Z})} = 0$ as well. However $V_n$ is generated by $V_a \subseteq V_n^{GL_{n-a}(\mathbb{Z})}$, so every irreducible component of $V_n$ (and therefore of $W_n / \widetilde{W}_n$) must contain a $GL_{n-a}(\mathbb{Z})$-invariant vector. So $W_n / \widetilde{W}_n = 0$, as desired.$\square$

\medskip

As a quick corollary, we obtain the following structural result about pointwise finite dimensional $VIC(Z)$-modules:

\medskip

\begin{corollary}\label{inject}
If $V$ is a finitely generated pointwise finite dimensional $VIC(\mathbb{Z})$-module then for $n >> 0$ the map $V_n \to V_{n+1}$ coming from the standard inclusion of $\mathbb{Z}^n$ into $\mathbb{Z}^{n+1}$ is injective.

\end{corollary}

\noindent \textbf{Proof:} The kernels of these maps form a $VIC(\mathbb{Z})$-submodule.  By Noetherianity this must be finitely generated, but all the maps are zero so these kernels must be zero in sufficiently large degree. $\square$

\medskip

The point at which these maps start being injective is referred to as the injectivity degree.

\subsection{Stability of depth for finite type representations}

The main result of this section is the following proposition, a structure theorem for $VIC(\mathbb{Z})$-modules with the additional assumption that the representations are pointwise of finite type.

\begin{proposition} \label{fintypevic}
Suppose $V$ is a finitely generated $VIC(\mathbb{Z})$-module which is pointwise finite dimensional and of finite type then there exists an $\ell$ such that each $V_n$ has depth dividing $\ell$ and $V$ has the structure of a $VIC^\pm(\mathbb{Z}/\ell \mathbb{Z})$-module.

\end{proposition}

The main idea of the proof will be to propagate the property that the level $\ell$ congruence subgroup acts trivially from one representation to the next.  In order to do so we'll need to introduce a bit of notation for working with a few adjacent terms in a $VIC$-module. Suppose we have three representations  $$V_{n-1} \to V_n \to V_{n+1}$$ such that:
\begin{enumerate}

\item Each $V_k$ is a representation of $GL_k(\mathbb{Z})$.

\item The maps are $GL_{n-1}(\mathbb{Z})$- and $GL_{n}(\mathbb{Z})$-equivariant respectively.

\item $V_{n}$ is generated by image of $V_{n-1}$ as a $GL_{n}(\mathbb{Z})$ representation, and $V_{n+1}$ is generated by image of $V_n$ as a $GL_{n+1}(\mathbb{Z})$ representation.

 \item The images of $V_{n-1}$ in $V_{n}$, and of $V_n$ in $V_{n+1}$ are fixed by the complementary copies of $GL_1(\mathbb{Z})$.

 \item The image of $V_{n-1}$ in $V_{n+1}$ is fixed by the complementary copy of $GL_2(\mathbb{Z})$.

\end{enumerate}
Call three such representations $V_{n-1}$, $V_n$, $V_{n+1}$ a \emph{weak $VIC(\mathbb{Z})$-triple.}  If additionally the maps are assumed to be injective we'll say they form a \emph{$VIC(\mathbb{Z})$-triple.} The point being that for a finitely generated pointwise finite dimensional $VIC(\mathbb{Z})$-module,  after the degree of generation every three adjacent terms form a weak $VIC(\mathbb{Z})$-triple, and after the both the degree of generation and the injectivity degree they form a $VIC(\mathbb{Z})$-triple.

 We'll use this notion a bit flexibly, we can change $\mathbb{Z}$ to another ring $R$, or we can replace $GL_n$ by $SL^\mathfrak{U}_n$ to obtain definitions for $VIC(R)$- and $VIC^\mathfrak{U}(R)$-triples as well.
 
 \medskip
 
We'll now state the key lemma which will let us propagate depth bounds down a $VIC(\mathbb{Z})$-module.
 
 \begin{lemma}\label{depthprop}  Suppose $$V_{n-1} \to V_n \to V_{n+1}$$ is a weak $VIC(\mathbb{Z})$-triple of representations such that the actions of $GL_{n-1}(\mathbb{Z})$ and $GL_{n}(\mathbb{Z})$ factor through $SL^\pm_{n-1}(\mathbb{Z}/\ell \mathbb{Z})$ and $SL^\pm_{n}(\mathbb{Z}/\ell \mathbb{Z})$ respectively.  Then the action of $GL_{n+1}(\mathbb{Z})$ on $V_{n+1}$ factors through $SL^\pm_{n+1}(\mathbb{Z}/\ell \mathbb{Z})$.
 \end{lemma}
 
 \noindent \textbf{Proof:} Since $V_{n+1}$ is generated by $V_{n-1}$ it suffices to check that the image of $V_{n-1}$ is fixed by the congruence subgroup $\Gamma_{n+1}(\ell)$.  To see this, note that it is fixed by $\Gamma_n(\ell)$ as it factors through $V_n$, and that it is also fixed by the complementary copy of $SL^\pm_2(\mathbb{Z})$ by the weak $VIC$-triple assumption, therefore it is fixed by the subgroup they generate. So it is enough to see that this subgroup generated by $\Gamma_n(\ell)$ and the complementary copy of $SL^\pm_2(\mathbb{Z})$ contains $\Gamma_{n+1} (\ell)$. 
 
  To see this, first note that the elementary matrices $E_{ij}^\ell$ (i.e. matrices with 1's on the diagonal, a single $\ell$ in position $(i,j)$, and zeroes everywhere else) lie in this subgroup - most of them lie in $\Gamma_n(\ell)$ already and those that aren't can be written as a commutator of one that is with $E_{n,n+1}$ or $E_{n,n+1}$ (which are in the complementary copy of $SL_2(\mathbb{Z})$).  In general these matrices $E_{ij}^\ell$ do not generate $\Gamma_{n+1} (\ell)$, but Theorem 7.5 of \cite{BMS} tells us that $\Gamma_{n+1} (\ell)$ is generated by the $E_{ij}^\ell$'s along with $\Gamma_{n} (\ell)$ (or even just $\Gamma_{2} (\ell)$). $\square$
 
\medskip 
 
\noindent \textbf{Proof of Proposition \ref{fintypevic}:} Suppose $V$ is pointwise finite dimensional and finite type, and is generated in degrees $\le d$. Let $\ell' = \text{l.c.m.}(\text{depth}(V_d), \text{depth}(V_{d+1}))$, for every $n >d $ each triple $$V_{n-1} \to V_n \to V_{n+1}$$ forms a weak $VIC(\mathbb{Z})$-triple.  Hence by Lemma \ref{depthprop} and induction we see that $\text{depth}(V_n)$ divides $\ell'$ for all $n >d$, and the action each $V_n$  factors through $SL^\pm_{n}(\mathbb{Z}/\ell' \mathbb{Z})$ for for $n>d$.

If we replace each $V_n$ with $n \le d$ by zero, Corollary \ref{pointwise} tells us that the resulting $VIC(\mathbb{Z})$-module is actually a $VIC(\mathbb{Z}/\ell' \mathbb{Z})^\pm$-module, which is almost what we want. To recover the proposition for all of $V$ we just need to increase $\ell'$ so that the first few terms before the degree of generation are accounted for.  In particular, we can take $\ell = \text{l.c.m}(\text{depth}(V_3), \text{depth}(V_{4}), \dots, \text{depth}(V_{d-1}), \ell' )$ and then again Corollary \ref{pointwise} tells us that $V$ has the structure of a $VIC^\pm(\mathbb{Z}/\ell \mathbb{Z})$-module. $\square$

\medskip

\noindent \textbf{Remark:} We'll also note that the same proof holds for a local version: if we replace $\mathbb{Z}$ by $\mathbb{Z}_p$ and look at pointwise smooth and finite dimensional $VIC(\mathbb{Z}_p)$-modules, then now the $p$-depth stabilizes and each such module is a pullback of a $VIC(\mathbb{Z}/p^k \mathbb{Z})$-module for some $k$. 

\medskip

\subsection{Extending from $SL^\mathfrak{U}_n$ to $GL_n$}
We have this mildly annoying issue that sometimes we have actions of $SL_n$, sometimes of $SL^\mathfrak{U}_n$, and sometimes of $GL_n$ (often over different rings).  One way to avoid this would be to just always restrict to $SL_n$ since it always acts, however the problem with doing that is that it is often easier to work with representations of general linear groups rather than of special linear groups.  Instead we will show that in the context representation stability there is a natural way to extend from $SL_n$ or $SL^\mathfrak{U}_n$ to $GL_n$ while satisfying the compatibility conditions between the representations.

 Over the complex numbers, every irreducible (algebraic) representation of $GL_n(\mathbb{C})$ remains irreducible when restricted to $SL_n(\mathbb{C})$, and every irreducible $SL_n(\mathbb{C})$ representation arises this way.  Different $GL_n(\mathbb{C})$ representations can give rise to the same $SL_n(\mathbb{C})$ representation, but two different extensions of an irreducible $SL_n(\mathbb{C})$ representation to $GL_n(\mathbb{C})$ differ by tensoring with a power of the determinant character. Over other rings the situation can be more complicated.

The main theorem of this section is that in the context of representation stability that not only can we extend each representation from $SL^\mathfrak{U}_n$ to $GL_n$, but we can do so for the different values of $n$ in a compatible way to get a stable sequence of general linear group representations.

\begin{theorem}\label{SLtoGL}

Suppose $V$ is a $VIC^\mathfrak{U}(R)$-module with degree of generation and injectivity degree $\le k$, there exists a $VIC(R)$-module $\hat{V}$ which restricts to $V$ in degrees $\ge k$.


\end{theorem}

\medskip

\noindent \textbf{Proof:} The condition that we are past the degrees of generation and injectivity ensures that each adjacent triple $$V_{n-1} \to V_{n} \to V_{n+1}$$ is a $VIC^\mathfrak{U}(R)$-triple for all $n > k$. The strategy will be to look at these triples and define actions of the general linear groups on each vector space  $V_n$ extending the actions of $SL_n^\mathfrak{U}(R)$, and to then check that new actions themselves form $VIC(R)$-triples. 

As a start note that we do automatically have an action of $GL_n(R)$ on $V_{n+1}$ coming from the embedding of $GL_n(R)$ in $SL_{n+1}(R)$ sending a matrix $A \in GL_n(R)$ to

 
 \[
 \begin{bmatrix}
   A      & 0 \\
    0       & \text{det}(A)^{-1} \\
 
\end{bmatrix}
\]

We claim that this action of $GL_n(R)$ on $V_{n+1}$ preserves the copy of $V_n$ in $V_{n+1}$ (we know the maps are injective since $n>k$ by assumption), and therefore defines an action of $GL_n(R)$ on $V_{n}$.  To see this, recall that our assumption of being a $VIC^\mathfrak{U}(R)$-triple gives us that $V_n$ is generated as an $SL_n^\mathfrak{U}(R)$ representation by $V_{n-1}$ which explicitly means any $\vec{v} \in V_n$ can be written as $\vec{v} =  \sum_i g_i \vec{v}_i$ with $g_i \in SL_n(R)$ and $\vec{v}_i \in V_{n-1}$.  Moreover, we know that the copy of $V_{n-1}$ in $V_{n+1}$ is fixed by all matrices of the form

 \[
 \begin{bmatrix}
   \text{Id}_{n-1}     & 0 & 0\\
    0       & t & 0 \\
    
    0 & 0 & t^{-1} \\
 
\end{bmatrix}
\]
Now note that any element of this copy of $GL_n(R)$ in $SL_{n+1}(R)$ can be written as a product 

 \[
 \begin{bmatrix}
   A      & 0 \\
    0       & \text{det}(A)^{-1} \\
\end{bmatrix} = 
 \begin{bmatrix}
   \hat{A}      & 0 \\
    0       & 1 \\
\end{bmatrix} \cdot
 \begin{bmatrix}
   \text{Id}_{n-1}     & 0 & 0\\
    0       & \text{det}(A) & 0 \\
    
    0 & 0 & \text{det}(A)^{-1} \\
 
\end{bmatrix}
\]
where $\hat{A} \in SL_n(R)$ is obtained from $A$ by dividing the last column of $A$ by $\text{det}(A)$ (Note: these are $n \times n$ matrices so this multiplication makes sense, despite how it looks).  Moreover, for any $g_i \in SL_n(R)$ we have that 

\[
 \begin{bmatrix}
   \text{Id}_{n-1}     & 0 & 0\\
    0       & \text{det}(A) & 0 \\
    
    0 & 0 & \text{det}(A)^{-1} \\
 
\end{bmatrix}\cdot
 \begin{bmatrix}
   g_i      & 0 \\
    0       & 1 \\
\end{bmatrix}=
 \begin{bmatrix}
   \tilde{g_i}      & 0 \\
    0       & 1 \\
\end{bmatrix} \cdot
 \begin{bmatrix}
   \text{Id}_{n-1}     & 0 & 0\\
    0       & \text{det}(A) & 0 \\
    
    0 & 0 & \text{det}(A)^{-1} \\
 
\end{bmatrix}
\]
where $\tilde{g_i}$ is obtained from $g_i$ by dividing the last column by $\text{det}(A)$ and then multiplying the last row by $\text{det}(A)$.  In particular for any $\vec{v} =  \sum_i g_i \vec{v}_i$ and $A \in GL_n(R)$ we can rewrite

 \[
 \begin{bmatrix}
   A      & 0 \\
    0       & \text{det}(A)^{-1} \\
\end{bmatrix} \vec{v}
\]
as

 \[
 \sum_i \begin{bmatrix}
   \hat{A}      & 0 \\
    0       & 1 \\
\end{bmatrix} \cdot
 \begin{bmatrix}
   \tilde{g_i}      & 0 \\
    0       & 1 \\
\end{bmatrix} \cdot
\begin{bmatrix}
   \text{Id}_{n-1}     & 0 & 0\\
    0       & \text{det}(A) & 0 \\
    
    0 & 0 & \text{det}(A)^{-1} \\
 
\end{bmatrix} \vec{v}_i = \sum_i
\begin{bmatrix}
   \hat{A} \tilde{g}_i     & 0 \\
    0       & 1 \\
\end{bmatrix} \vec{v}_i
\]
Which is clearly in the space generated under the action of $SL_n(R)$ by $V_{n-1}$, and we know that space is $V_n$ as they are part of a $VIC^\mathfrak{U}(R)$-triple. So indeed the copy of $V_n$ in $V_{n+1}$ is preserved by this copy of $GL_n(R)$ in $SL_{n+1}$, as desired.

Now that we have these actions of $GL_n(R)$ defined pointwise for each $n$ we still need to check that they are compatible with one another and form $VIC(R)$-triples.  First, we need to show that the map $V_{n-1} \to V_{n}$ is $GL_{n-1}(R)$-equivariant \emph{for the standard embedding of $GL_{n-1}(R)$ in $GL_{n}(R)$}, since right now it is defined in terms of the action of a different copy of $GL_{n-1}(R)$ in $GL_n(R)$.  To see this, we need to look inside $V_{n+1}$. 

We have two copies of $GL_{n-1}(R)$ in $SL_{n+1}(R)$: The first is given by the embedding of $GL_{n-1}(R)$ into $SL_n(R)$ from above followed by the standard embedding of $SL_n(R)$ into $SL_{n+1}(R)$.  Explicitly this sends a matrix $A$ to

 \[
 \begin{bmatrix}
   A      & 0 & 0\\
    0       & \text{det}(A)^{-1} & 0 \\
    
    0 & 0 & 1\\
 
\end{bmatrix}
\]

The second copy of $GL_{n-1}(R)$ in $SL_{n+1}(R)$ we care about is the one obtained by the standard embedding of $GL_{n-1}(R)$ in $GL_{n}(R)$ followed by the embedding of $GL_{n}(R)$ in $SL_{n+1}(R)$ similar to what we had before.  Explicitly here a matrix $A$ gets sent to

 \[
 \begin{bmatrix}
   A      & 0 & 0\\
    0       & 1 & 0 \\
    
    0 & 0 & \text{det}(A)^{-1} \\
 
\end{bmatrix}
\]

Now note that these matrices may act differently on $V_n$ and $V_{n+1}$ (and most likely do).  However in order for our map $V_{n-1} \to V_n$ to be $GL_{n-1}(R)$ equivariant we just need them to act the same on the copy of $V_{n-1}$ inside $V_{n+1}$. This is true since these two matrices are conjugate via the matrix 

 \[
 \begin{bmatrix}
   \text{Id}_{n-1}     & 0 & 0\\
    0       & 0 & -1 \\
    
    0 & 1 & 0 \\
 
\end{bmatrix}
\]
which acts trivially on $V_{n-1}$ because of our central stability assumption.

We also see by this calculation that the matrices

  \[
 \begin{bmatrix}
   A      & 0 \\
    0       & \text{det}(A)^{-1}
\end{bmatrix} \text{ and }
 \begin{bmatrix}
   A      & 0 \\
    0       & 1 \\
 
\end{bmatrix}
\]
acting on $V_n$ have the same action when restricted to $V_{n-1}$.  In particular multiplying one by the inverse of the other tells us that any matrix 

 \[
 \begin{bmatrix}
   \text{Id}     & 0 \\
    0       & t \\
 
\end{bmatrix}
\]
acts trivially on $V_{n-1}$, so indeed $V_{n-1}$ is fixed by the complementary copy of $GL_1(R)$ in $GL_n(R)$.

The facts that the maps are injective, and that the image of $V_{n-1}$ generates $V_{n}$ as a $GL_n(R)$ representation for each $n$ follows for free from the fact that it already generated for $SL(R)$ representations. So in order to prove each adjacent triple is actually a $VIC$-triple all that is left is to verify is that $V_{n-1}$ is preserved by the complementary copy of $GL_2(R)$ inside $GL_{n+1}(R)$.

We know that $V_{n-1}$ is preserved by $SL_2(R)$ and it's also preserved by the matrix

 \[
 \begin{bmatrix}
   \text{Id}_{n-1}     & 0 & 0\\
    0       & 1 & 0 \\
    
    0 & 0 & t \\
 
\end{bmatrix}
\]
for any $t$ since $V_{n-1}$ is inside $V_n$, which is preserved by this by the previous calculation.  The complementary $GL_2(R)$ is generated by $SL_2(R)$ and matrices of this form, so indeed the full $GL_2(R)$ acts trivially on this space. $\square$

\medskip

\noindent \textbf{Remark:} This theorem is analogous to a result of Wilson which says that every $FI_D$-module is, in large degrees, a restriction of an $FI_{BC}$-module (Prop. 3.30 in \cite{Wil}). However the proof is quite different.

\medskip

One of the main reasons for doing this is the following corollary:

\begin{corollary}\label{algextend} \ 

\begin{enumerate}

\item A finitely generated $VIC(\mathbb{Z})$-module that is pointwise finite dimensional and finite type has the structure of a finitely generated $VIC(\mathbb{Z} / \ell \mathbb{Z})$-module in large enough degrees.

\item A finitely generated $VIC(\mathbb{Z})$-module that is pointwise finite dimensional and algebraic is the restriction of a finitely generated pointwise algebraic $VIC(\mathbb{C})$-module in large enough degrees.

\item More generally, the induced $VIC^\pm(\mathbb{C})$- and $VIC^\pm(\mathbb{Z}_p)$-module structures on a finitely generated pointwise finite dimensional $VIC(\mathbb{Z})$-module from Proposition \ref{weaksr} can be extended to $VIC(\mathbb{C})$- and $VIC(\mathbb{Z}_p)$-module structures in large degrees.

\end{enumerate}
\end{corollary}

\noindent \textbf{Proof:} For parts 1 and 2,  Proposition \ref{fintypevic} and Corollary \ref{pointwise} tell us that these have the structure of a $VIC^\pm(\mathbb{Z} / \ell \mathbb{Z})$ and a pointwise algebraic $VIC^\pm(\mathbb{C})$-module respectively. Corollary \ref{inject} tells us that finitely generated pointwise finite dimensional $VIC(\mathbb{Z})$-modules have injective maps in large enough degrees, so Theorem \ref{SLtoGL} applies and we can therefore extend them to $VIC(\mathbb{Z} / \ell \mathbb{Z})$ and $VIC(\mathbb{C})$ respectively.

For part 3, we'll once again note that in general the induced $VIC^\pm(\mathbb{C})$ and $VIC^\pm(\mathbb{Z}_p)$-module structures will not be finitely generated, so we can't apply Theorem \ref{SLtoGL} directly to either of them individually. We can however apply it to them simultaneously by using superrigidity to interpret our $VIC(\mathbb{Z})$-module as a finitely generated $VIC^\pm(R)$-module for $R = \mathbb{C} \times \prod_p \mathbb{Z}_p$. Again  Corollary \ref{inject} gives us injectivity in large degrees so Theorem \ref{SLtoGL} then lets us extend from $SL^\pm(\mathbb{C}) \times \prod_p SL^\pm(\mathbb{Z}_p)$ to $GL(\mathbb{C}) \times \prod_p GL(\mathbb{Z}_p)$ in large degrees, and restricting to one of the factors gives the desired result.  $\square$

\medskip

\noindent \textbf{Remark:} A special case of part 2 of the above corollary appears as Patzt's Theorem 5.3 in \cite{Patzt1}. 

\medskip

\subsection{Algebraic parts of $VIC(\mathbb{Z})$-modules}

Patzt has previously shown that a finitely generated $VIC(\mathbb{C})$-module that is pointwise algebraic exhibits representation stability in the sense that under the standard labeling of irreducibles by bipartitions (described below) the multiplicity of an irreducible corresponding to a fixed bipartition eventually stabilizes, and only finitely many bipartitions occur (see \cite{Patzt1} Theorem A).

We will be drawing on many of Patzt's ideas in order to get a hold on the algebraic parts of the $GL_n(\mathbb{Z})$-representations that arise in a finitely generated pointwise finite dimensional $VIC(\mathbb{Z})$-module. Before we prove the main theorem we will review a bit about algebraic representations of $GL_n(\mathbb{C})$, and the ingredients of Patzt's proof that we will be drawing on.

\subsubsection{Algebraic representations of $GL_n(\mathbb{C})$}

First a bit of background on algebraic representations of $GL_n(\mathbb{C})$.  Highest weight theory for a reductive algebraic group gives a natural labeling of irreducible finite dimensional algebraic representations in terms of dominant integral weights. In the case of $GL_n(\mathbb{C})$ these are just $n$-tuples $(a_1,a_2,\dots, a_n)$ of integers which are weakly decreasing $a_1 \ge a_2 \ge \dots \ge a_n$.

For a dominant integral weight $\omega$ let $W^\omega$ denote the corresponding irreducible representation. For example, under this labeling the standard $n$ dimensional representation is $W^{(1,0,0\dots 0)}$, its dual is $W^{(0,0,\dots, 0, -1)}$, and the conjugation action on $\mathfrak{sl}_n$ is $W^{(1,0\dots,0, -1)}$.

In the settings of representation stability, Deligne categories, and supergroup representations it is often convenient to use a slightly different parameterization in terms of pairs of partitions $(\lambda^+, \lambda^-)$, which we will call bipartitions.  If $\lambda^+ = \lambda^+_1 \ge \lambda^+_2 \ge... \ge \lambda^+_r$ and $\lambda^- = \lambda^-_1 \ge \lambda^-_2 \ge... \ge \lambda^-_s$, then for $n \ge r+s$ define $$\omega(\lambda^+,\lambda^-):= (\lambda^+_1, \lambda^+_2, \dots, \lambda^+_r, 0 ,\dots, 0, -\lambda^-_s, \dots, -\lambda^-_2, -\lambda^-_1) $$ and we'll call $V_n(\lambda^+, \lambda^-) = W^{\omega(\lambda^+,\lambda^-)}$.  Now, for all $n$ the defining $n$-dimensional representation is called $V_n(\square, \emptyset)$, its dual is $V_n(\emptyset, \square)$, the $\mathfrak{sl}_n$ representation is $V_n(\square,\square)$, and the trivial representation is $V_n(\emptyset, \emptyset)$.

The irreducible representations of the form $V_n(\lambda, \emptyset)$ play a special role in this theory, they are the so called polynomial representations of $GL_n(\mathbb{C})$. Alternatively these are those $W^\omega$ where $\omega = (a_1,a_2,\dots, a_n)$ with $a_n \ge 0$. They are exactly those representations such that the map $GL_n(\mathbb{C}) \to GL_N(\mathbb{C})$ is given by polynomial functions of the matrix entries (as opposed to rational functions).

 The irreducible polynomial representation $V(\lambda, \emptyset)$ is obtained by applying the Schur functor $S^\lambda$ to the standard $n$-dimensional representation $V(\square, \emptyset)$.  The characters of these representations are given by symmetric polynomials, with $V(\lambda, \emptyset)$ corresponding to the Schur polynomial $s_\lambda$.

One key fact is that if $\omega = (a_1,a_2,\dots, a_n)$ is a dominant integral weight, then $W^\omega \otimes \text{det} = W^{\omega + \vec{1}}$ where $\text{det}$ denotes the determinant representation and $\omega + \vec{1} = (a_1+1,a_2+1,\dots, a_n+1)$. In particular, one can always tensor an algebraic representation with some power of the determinant representation to obtain a polynomial representation.

\subsubsection{Branching of algebraic representations}

For polynomial representations the branching rules for restricting from $GL_{n+m}(\mathbb{C})$ to $GL_n(\mathbb{C}) \times GL_m(\mathbb{C})$ are well known: they are just given by Littlewood-Richardson coefficients.  Explicitly, the multiplicity of $V_n(\lambda,\emptyset) \otimes V_m(\mu, \emptyset)$ in the restriction of $V(\nu,\emptyset)$ is $c^\nu_{\lambda,\mu}$ provided $n$ and $m$ are larger than the lengths of $\lambda$ and $\mu$ respectively (if not, those partitions don't correspond to weights).

Using this one can easily deduce branching rules for algebraic representations by tensoring with an appropriate power of the determinant to have the representation become polynomial, and then keeping track of the degree shift. We won't go into detail into the combinatorics of the modified Littlewood-Richardson rule for algebraic representations, but to give a flavor we will state the simplest case -- the Pieri rule for algebraic representations.  Let for a partition $\lambda$ let $\ell(\lambda)$ denote its length (i.e. the number of parts).  
 
 \begin{lemma}[Pieri rule for algebraic representations] Suppose $\ell(\lambda^+) + \ell(\lambda^-) < n$ then we have
 
$$\text{Res}^{GL_n(\mathbb{C})}_{GL_{n-1}(\mathbb{C})\times \mathbb{C}^{\times}}V_n(\lambda^+, \lambda^-) \cong \bigoplus_{(\mu^+,\mu^-) \in \text{hs}(\lambda^+,\lambda^-)} V_{n-1}(\mu^+,\mu^-)\otimes 
\mathbb{C}_{(|\lambda^+|-|\mu^+|-|\lambda^-|+|\mu^-|)}$$
 where $\text{hs}(\lambda^+,\lambda^-)$ denotes the set of pairs of partitions obtained from $(\lambda^+,\lambda^-)$ by removing a horizontal strip from each side, and $\mathbb{C}_k = W^{(k)}$ is the character of $\mathbb{C}^{\times} = GL_1(\mathbb{C})$ where $t$ acts by $t^k$.
 
\end{lemma}
 
 
 In particular, consider the containment partial order on bipartitions where $(\mu^+, \mu^-) \le (\lambda^+, \lambda^-)$ if and only if $\mu^+$ is contained in $\lambda^+$ and $\mu^-$ is contained in $\lambda^-$ (by containment we mean containment of the corresponding Young diagrams).  It's clear from the Pieri rule that if $\ell(\lambda^+) + \ell(\lambda^-) \le m$ then the restriction of $V_n(\lambda^+, \lambda^-)$ from $GL_n(\mathbb{C})$ to $GL_m(\mathbb{C})$ only has terms $V_m(\mu^+, \mu^-)$ with  $(\mu^+, \mu^-) \le (\lambda^+, \lambda^-)$ in the containment order.
 
 \medskip
 
We will only need a couple key facts about these branching rules, which were also important to Patzt's analysis of pointwise algebraic $VIC(\mathbb{C})$-modules.  The first is Patzt's Lemma 3.16 in \cite{Patzt1}, which just uses the Pieri rule and gives control over the length of partitions appearing in representations with $GL_{n-m}(\mathbb{C})$-invariant vectors:
 
 \begin{lemma}\label{branchingcor}Suppose $\ell(\lambda^+) + \ell(\lambda^-) < n$

 \begin{enumerate}  
 
  \item Every term $V_{n-1}(\mu^+,\mu^-)\boxtimes \mathbb{C}_{(|\lambda^+|-|\mu^+|-|\lambda^-|+|\mu^-|)}$ that appears the Pieri rule restriction of $V_n(\lambda^+,\lambda^-)$ has $\ell(\lambda^+) - 1 \le \ell(\mu^+) \le \ell(\lambda^+)$ and $\ell(\lambda^-) - 1 \le \ell(\mu^-) \le \ell(\lambda^-)$.
  
  \item In particular, if the restriction of $V(\lambda^+, \lambda^-)$ to $GL_{n-1}(\mathbb{C})$ has any trivial factors $V(\emptyset, \emptyset)$ then necessarily $\ell(\lambda^+) + \ell(\lambda^-) \le 2$.  Iterating this, it follows that if the restriction of $V_n(\lambda^+,\lambda^-)$ to $GL_{n-m}(\mathbb{C})$ has any trivial factors then $\ell(\lambda^+) + \ell(\lambda^-) \le 2m$.
  
\end{enumerate}  
  
\end{lemma}
 
 The second fact we need is another special case of the modified Littlewood-Richardson rule and is Patzt's Proposition 3.9 in \cite{Patzt1} and is in a sense a converse to the previous lemma. This says that provided the lengths of the partitions involved are small enough for the corresponding representations to be defined, irreducible algebraic representations with the same label actually appear in the appropriate invariant subspaces of one another.
 
 \begin{lemma} \label{branching2} Suppose $\ell(\lambda^+) + \ell(\lambda^-) \le m < n$. The restriction from $GL_{n}(\mathbb{C})$  to $GL_{m}(\mathbb{C}) \times GL_{n-m}(\mathbb{C})$ of $V_n(\lambda^+, \lambda^-))$ contains a unique copy of $V_m(\lambda^+, \lambda^-) \otimes V_{n-m}(\emptyset, \emptyset)$.  In other words, $V_n(\lambda^+, \lambda^-))^{GL_{n-m}(\mathbb{C})}$ contains a unique copy of $V_m(\lambda^+, \lambda^-)$ as a representation of $GL_{m}(\mathbb{C})$.

\end{lemma}

\subsubsection{A filtration of $VIC(\mathbb{Z})$-modules by algebraic part}

If $M$ is a representation of $GL_n(\mathbb{C})$ let $M^{\ge (\lambda^+, \lambda^-)}$ denote the subspace of $M$ generated by all irreducible components isomorphic to $V^{(\mu^+ ,\mu^-)}$ where $\mu^+$ contains $\lambda^+$ and $\mu^-$ contains $\lambda^-$.  An important observation of Patzt, and a consequence of the Pieri rule, is that if $V$ is a pointwise algebraic $VIC(\mathbb{C})$-module then the subspaces $V_n^{\ge (\lambda^+, \lambda^-)}$ form a $VIC(\mathbb{C})$-submodule.

We'd now like to extend this to our setting of $VIC(\mathbb{Z})$-modules. A priori superrigidity for a single $GL_n(\mathbb{Z})$ only tells us we can virtually extend a representation to $SL_n^\pm(\mathbb{C})$, however in light of Corollary \ref{algextend} we know in large enough degrees this extends naturally to $GL_n(\mathbb{C})$. In particular, the the subspaces $V_n^{\ge (\lambda^+, \lambda^-)}$ are still defined in large degrees, and still form a $VIC(\mathbb{Z})$-submodule. 

Rephrasing Patzt's argument and adapting it a bit for our setting, we can then isolate the components of $V_n$ where the algebraic part is $V_n(\lambda^+,\lambda^-)$ for a fixed bipartition $(\lambda^+, \lambda^-)$ by taking a $VIC(\mathbb{Z})$-subquotient. Explicitly this gives us the following lemma:

\begin{lemma}\label{algisotypic}

A finitely generated pointwise finite dimensional $VIC(\mathbb{Z})$-module $V$ contains submodules $V^{\ge (\lambda^+, \lambda^-)}$ and $V^{> (\lambda^+, \lambda^-)}$ such that for all sufficiently large $n$: 

\begin{enumerate}

\item $(V^{ \ge (\lambda^+, \lambda^-)} / V^{> (\lambda^+, \lambda^-)})_n \cong V_n(\lambda^+, \lambda^-) \otimes W_n$ where $W_n$ is a smooth representation of $\prod_p GL_n(\mathbb{Z}_p)$.

\item  $(V/V^{ \ge (\lambda^+, \lambda^-)})_n$ does not contain any copies of $V_n(\lambda^+, \lambda^-)$ as $GL_n(\mathbb{C})$-representations

\item $V^{> (\lambda^+, \lambda^-)}_n$ contains only irreducible representations with algebraic parts $V_n(\mu^+, \mu^-)$ where $(\mu^+, \mu^-)$ is striuctly larger than $(\lambda^+, \lambda^-)$ in the containment order on bipartitions.

\end{enumerate}

\end{lemma}

\medskip

In Patzt's case where the representations are pointwise algebraic, the actions on the $W_n$'s are all automatically trivial. In that case Patzt then proves that the dimensions of these spaces $W_n$ stabilize in $n$, and moreover that only finitely many bipartitions $(\lambda^+,\lambda^-)$ appear with nonzero multiplicity.

\medskip

\subsection{The structure of $VIC(\mathbb{Z})$-modules}

Now we'd like to give a description of a finitely generated $VIC(\mathbb{Z})$-module as a sequence of $GL_n(\mathbb{Z})$ representations.  

\begin{theorem} \label{main}
Let $V$ be a finitely generated pointwise finite dimensional $VIC(\mathbb{Z})$-module over $\mathbb{C}$,  $V$ admits a finite filtration $$V \supseteq V^1 \supseteq \dots \supseteq V^k = 0$$ by $VIC(\mathbb{Z})$-submodules such that $V / V^1$ is a torsion $VIC(\mathbb{Z})$-module (meaning $(V / V^1)_n = 0$ for all $n \gg 0$) and for each other successive subquotient $V^i / V^{i+1}$ there exists a bipartition $(\lambda^+_i, \lambda^-_i)$ such that $$(V^i / V^{i+1})_n = V_n(\lambda^+_i, \lambda^- _i) \otimes M^i_n$$ 
where the collection of $M^i_n$'s form a finitely generated $VIC(\mathbb{Z} / \ell_i \mathbb{Z})$-module for some integer $\ell_i$.

\end{theorem}

\noindent \textbf{Proof:} We've broken down the proof into 4 steps:

\medskip

\noindent \textbf{Step 1: Throw out the low degree terms.} Corollary \ref{algextend} tells us that any finitely generated pointwise finite dimensional $VIC(\mathbb{Z})$-module over $\mathbb{C}$ has the structures of being $VIC(\mathbb{C})$ and $VIC(\mathbb{Z}_p)$ modules (for all $p$) in large enough degrees.  We'll take $V^1$ to just be $V$ in those degrees for which the corollary applies, and zero in low degrees. So from now on we may always assume $n$ is large enough so that each $V_n$ has induced actions of $GL_n(\mathbb{C})$ and $GL_n(\mathbb{Z}_p)$.

\medskip

\noindent \textbf{Step 2: Define the filtration.} If we fix a linear refinement of the containment order then Lemma \ref{algisotypic} gives us a descending filtration by $VIC(\mathbb{Z})$-submodules, where each successive quotient $V^i / V^{i+1}$ is pointwise given by $$(V^i / V^{i+1})_n = V_n(\lambda^+_i, \lambda^- _i) \otimes M^i_n$$ where $(\lambda^+_i, \lambda^- _i)$ is the $(i+1)$th bipartition in the linear refinement, and $M^i_n$ is a (possibly zero) finite type representation which we know extends to $GL_n(\mathbb{Z} / \ell \mathbb{Z})$ for some $\ell$ (a priori depending on $i$ and $n$).  Our final filtration will be this but reindexed to only include those bipartitions that actually appear.

\medskip

In order to finish the proof we still need to check that only finitely many bipartitions appear with nonzero multiplicities so that this is a finite filtration, and to check that for each bipartition $(\lambda^+_i, \lambda^- _i)$ the multiplicity spaces $M^i_n$ form a finitely generated $VIC(\mathbb{Z} / \ell_i \mathbb{Z})$-module for some $\ell_i$.

\medskip

\noindent \textbf{Step 3: Prove the filtration is finite.} We need to show that only finitely many bipartitions occur with non-zero multiplicities. Patzt proved this in the case of pointwise algebraic $VIC(\mathbb{C})$-modules, and we will follow his approach.  First note that if $V$ is generated in degrees $\le d$ then every $V_n$ is generated by its space of $GL_{n-d}(\mathbb{C})$ invariants, in particular this means that all bipartitions appearing in $V_n$ have $\ell(\lambda^+) + \ell(\lambda^-) \le 2d$ by Lemma \ref{branchingcor}.

If $V_n(\lambda^+, \lambda^-)$ appears as the algebraic part of a component of $V_n$ then Lemma \ref{branching2} tells us $V_{2d}(\lambda^+, \lambda^-)$ occurs inside $(V_n)^{GL_{n-2d}(\mathbb{Z})}$. By semisimplicity, we may identify the invariants $(V_n)^{GL_{n-2d}(\mathbb{Z})}$ with the covariants $(V_n)_{GL_{n-2d}(\mathbb{Z})}$.  But now we know by Corollary \ref{Phia} that these spaces of covariants stabilize as $GL_{2d}(\mathbb{Z})$ representations, and hence only finitely many bipartitions $(\lambda^+, \lambda^-)$ can appear.

\medskip

\noindent \textbf{Step 4: Prove the $M^i$ are finitely generated $VIC(\mathbb{Z})$-modules.}  So far we have finitely many subquotient $VIC(\mathbb{Z})$-modules $(V^i / V^{i+1})$, which pointwise look like $$(V^i / V^{i+1})_n = V_n(\lambda^+_i, \lambda^- _i) \otimes M^i_n$$ with $M^i_n$ of finite type, and by Noetherianity (Theorem \ref{Noetherianity}) we know that these are finitely generated. We now want to prove whenever we have a finitely generated $VIC(\mathbb{Z})$-module of this form that the multiplicity spaces $M^i_n$ themselves have the structure of a finitely generated $VIC(\mathbb{Z} / \ell_i \mathbb{Z})$-module for some $\ell_i$.

Note that we already have proven an important special case of this:  The finite type components of a finitely generated pointwise finite dimensional $VIC(\mathbb{Z})$-module (corresponding to $(\lambda^+_1, \lambda^- _1) = (\emptyset, \emptyset)$) form a $VIC(\mathbb{Z})$-quotient module since $(\emptyset, \emptyset)$ is the minimal element in the containment order. This is a finitely generated and pointwise finite type $VIC(\mathbb{Z})$-module therefore has the structure of a finitely generated $VIC(\mathbb{Z} / \ell_1 \mathbb{Z})$-module for some $\ell_1$ by part 1 of Corollary \ref{algextend}.

We will now reduce the general case to this special case. Take a finitely generated pointwise algebraic $VIC(\mathbb{Z})$-module $W^i$ such that $W^i_n = V_n(\lambda^-, \lambda^+)$ for $n \ge \ell(\lambda^+)+ \ell(\lambda^-)$ and consider the pointwise tensor product $W^i \otimes (V^i / V^{i+1})$.  Since $V_n(\lambda^-, \lambda^+)$ and $V_n(\lambda^+, \lambda^-)$ are dual this means that $V_n(\lambda^-, \lambda^+) \otimes V_n(\lambda^+, \lambda^-)$ contains a one dimensional invariant subspace (and of course the rest of the components are algebraic).  Therefore $M_n^i$ is the finite type component of $$(W^i \otimes (V^i / V^{i+1}))_n = V_n(\lambda^-, \lambda^+) \otimes V_n(\lambda^+, \lambda^-) \otimes M_n^i$$ but by Corollary \ref{tensor} we know that the pointwise tensor product has the structure of a finitely generated pointwise finite dimensional $VIC(\mathbb{Z})$-module, and is therefore by the special case above we see that the $M_n^i$ form a finitely generated $VIC(\mathbb{Z} / \ell_i \mathbb{Z})$-module for some $\ell_i$. $\square$

\medskip

The following corollary was conjectured by Putman and was one of the initial goals in starting this project, it is immediate from the above theorem:
\begin{corollary}
The sequence of representations coming from a finitely generated pointwise finite dimensional $VIC(\mathbb{Z})$-modules has bounded depth.
\end{corollary}

We can also now say a bit about the dimensions of these representations. Since the index of $GL_a(\mathbb{Z} /\ell \mathbb{Z}) \times GL_{n-a}(\mathbb{Z} /\ell \mathbb{Z})$ in $GL_n(\mathbb{Z} /\ell \mathbb{Z})$ grows exponentially in $n$ for fixed $a$ and $\ell$, and since the dimension of $V_n(\lambda^+, \lambda^-)$ is polynomial in $n$ we obtain the following:

\begin{corollary} \label{dimgrowth}
For any finitely generated $VIC(\mathbb{Z})$-module $V$ which is pointwise finite dimensional and finite type there is a constant $C$ such that $\text{dim}(V_n) \le C^n$ for all $n$.

\end{corollary}

\noindent \textbf{Remark:} We suspect that in fact the Hilbert series encoding these dimensions is rational with denominator a product of terms $(1 - jt)$ with $j \in \mathbb{N}$. Since $\text{dim}(V_n(\lambda^+, \lambda^-))$ is polynomial in $n$, this would follow from Theorem \ref{main} and an analogous statement for finitely generated $VIC(\mathbb{Z}/\ell \mathbb{Z})$-modules.  For the special class of $VI(\mathbb{Z}/\ell \mathbb{Z})$-modules this is known explicitly due to work of Sam and Snowden (\cite{SS1} Corollary 8.3.4).  It has been communicated to the author by Steven Sam that the same should be true for $VIC(\mathbb{Z}/\ell \mathbb{Z})$ by combining Lemma 2.15 in \cite{PS} with the methods of \cite{SS1}.

\medskip

We'd like to think of Theorem \ref{main} as a sort of multiplicity stability statement for pointwise finite dimensional $VIC(\mathbb{Z})$-modules. If the finite parts are all trivial we have the multiplicity stability statement of Patzt, and if all the representations factor through $GL_n(\mathbb{Z}/p\mathbb{Z})$ Gan and Watterlond have also shown a multiplicity stability statement (\cite{GW} Theorem 2).

 The reason that we can't literally state it as a multiplciity stability statement is that we don't know of a consistent naming system for the irreducible representations of $GL_n(\mathbb{Z}/\ell \mathbb{Z})$ for different $n$ (or really any natural labeling of these irreducible representations). We believe that such a naming system should exist and that finitely generated $VIC(\mathbb{Z} / \ell \mathbb{Z})$ over $\mathbb{C}$ should exhibit multiplicity stability, however for the time being we will settle for the following weaker statement.
 
 \begin{lemma}
 
 Suppose $V$ is a finitely generated $VIC(\mathbb{Z} / \ell \mathbb{Z})$ over $\mathbb{C}$. There exists an $N$ such that $V_n$ decomposes as a direct sum of at most $N$ irreducible $GL_n(\mathbb{Z} / \ell \mathbb{Z})$-representations for all $n > 2$.
 
 \end{lemma}

\noindent \textbf{Proof:}  If $V$ is generated in degrees $\le d$ then every irreducible component of each $V_n$ appears as a summand of $$\bigoplus_{i\le d}M(i)_n$$ where $M(i)$ denotes the free $VIC(\mathbb{Z} / \ell \mathbb{Z})$-module generated in degree $i$. This means that the dimension of $$Hom_{GL_n(\mathbb{Z} / \ell \mathbb{Z})}(\bigoplus_{i\le d}M(i)_n, V_n)$$ is a weighted sum of the multiplicities of the irreducible representations appearing in $V_n$, weighted by their corresponding multiplicities in $\bigoplus_{i\le d}M(i)_n$. In particular it is an upper bound for the number of components in a direct sum decomposition of $V_n$.  On the other hand, each $M(i)$ is self dual so this dimension is the same as the dimension of $$((\bigoplus_{i\le d}M(i)_n) \otimes V_n)^{GL_n(\mathbb{Z} / \ell \mathbb{Z})}$$ but these spaces $(\bigoplus_{i\le d}M(i)_n) \otimes V_n$ naturally form a finitely generated $VIC(\mathbb{Z} / \ell \mathbb{Z})$-module so we know that the dimension of the covariants, and hence the invariants since we are over $\mathbb{C}$ stabilizes. $\square$

\medskip

Combining this with Theorem \ref{main} immediately gives us an analogous statement for pointwise finite dimensional $VIC(\mathbb{Z})$-modules.

 \begin{corollary} \label{lengthbound}
 
 Suppose $V$ is a finitely generated pointwise finite dimensional $VIC(\mathbb{Z})$ over $\mathbb{C}$. There exists an $N$ such that $V_n$ decomposes as a direct sum of at most $N$ irreducible $GL_n(\mathbb{Z})$-representations for all $n > 2$.
 
 \end{corollary}

\section{Applications and extensions}\label{applications}

\subsection{Infinite rank superrigidity}\label{infrankform}

Let $GL_\infty(\mathbb{Z}) = \varinjlim GL_n(\mathbb{Z})$ be the group of one directionally infinite invertible matrices with integer entries which differ from the identity matrix in only finitely many positions. The inclusions $GL_k(\mathbb{Z}) \times GL_{n-k}(\mathbb{Z}) \hookrightarrow GL_n(\mathbb{Z})$ gives us an inclusion $GL_k(\mathbb{Z})\times GL_\infty(\mathbb{Z}) \hookrightarrow GL_\infty(\mathbb{Z})$.  We will denote the image of this second factor as $GL_{\infty - k}(\mathbb{Z})$.

Since by the maps are eventually injective by Corollary \ref{inject} the tail of a finitely generated pointwise finite dimensional $VIC(\mathbb{Z})$-module can be naturally thought of as a $GL_\infty(\mathbb{Z})$-representation by taking $V_\infty = \varinjlim V_n$. Since Theorem \ref{main} is ultimately about these tails, we may interpret it as a superrigidity theorem for a class of $GL_\infty(\mathbb{Z})$ representations.

Let $V$ be a representation of $GL_\infty(\mathbb{Z})$. We say that $V$ is \emph{admissible} if the following properties hold:

\begin{enumerate}

\item Every vector $\vec{v} \in V$ is fixed by $GL_{\infty - k}(\mathbb{Z})$ for some $k \in \mathbb{N}$.

\item The space $V^{GL_{\infty - k}(\mathbb{Z})}$ is finite dimensional for each $k \in \mathbb{N}$.

\item There exists a $k \in \mathbb{N}$ such that $V$ is generated by $V^{GL_{\infty - k}(\mathbb{Z})}$.  

\end{enumerate}

Any admissible $GL_\infty(\mathbb{Z})$-representation $V_\infty$ gives rise to a (torsion-free) $VIC(\mathbb{Z})$-module $V$, by taking $V_n = V_\infty^{GL_{\infty-n}(\mathbb{Z})}$. Similarly any $VIC(\mathbb{Z})$-module gives rise to an admissible $GL_\infty(\mathbb{Z})$ representation by taking $\varinjlim V_n$.  This defines an equivalence of categories between admissible $GL_\infty(\mathbb{Z})$-representations and the Serre quotient of the category of pointwise finite dimensional $VIC(\mathbb{Z})$-modules by the subcategory of torsion modules.

This gives us a dictionary for translating our results about pointwise finite dimensional $VIC(\mathbb{Z})$-modules into results about admissible $GL_\infty(\mathbb{Z})$ representations.   Sam and Snowden have given similar dictionaries for a number of other settings in representation stability, see \cite{SS2}. We'll now translate some of our results into this language. 

Under this correspondence, our Noetherianity result of Theorem \ref{Noetherianity} tells us the following:

\begin{theorem}

Any $GL_\infty(\mathbb{Z})$-subrepresentation of an admissible representation is also admissible, and therefore admissible $GL_\infty(\mathbb{Z})$-representations and homomorphisms between them form an abelian category.

\end{theorem}

We say that an admissible $GL_\infty(\mathbb{Z})$-representation is locally algebraic if any vector $\vec{v}$ generates an algebraic representation of $GL_n(\mathbb{Z})$ for all $n$, and that an admissible $GL_\infty(\mathbb{Z})$-representation is finite type if it factors through $GL_\infty(\mathbb{Z}/\ell \mathbb{Z})$ for some $\ell$.  Theorem \ref{main} may be rephrased as:

\begin{theorem}[Superrigidity for $GL_\infty(\mathbb{Z})$] \label{infrank}

Any admissible $GL_\infty(\mathbb{Z})$ representation has a finite filtration with subquotients of the form $V \otimes W$ where $V$ is locally algebraic and $W$ is finite type.

\end{theorem}

Finally our Corollary \ref{lengthbound} which bounded the number of factors in direct sum decompositions of each $V_n$ in a finitely generated pointwise finite dimensional $VIC(\mathbb{Z})$-module gives us the following corollary telling us that the category of admissible $GL_\infty(\mathbb{Z})$-representations is Krull-Schmidt.

\begin{corollary}

Any admissible $GL_\infty(\mathbb{Z})$ representation has finite length as a representation of $GL_\infty(\mathbb{Z})$.

\end{corollary}

It seems natural to ask if we can similarly write the the category of admissible $GL_\infty(\mathbb{Z})$-representations as the Deligne tensor product of the categories of admissible $GL_\infty(\mathbb{\mathbb{C}})$ and $GL_\infty(\mathbb{Z}_p)$ representations analogously to the formulation of superrigidity in Section \ref{catsup}. We suspect something like this should be true, but for the time being it seems out of reach.

\subsection{$VIC(\mathbb{Z})$-modules of small dimension growth}

We have shown in Corollary \ref{dimgrowth} that the dimension growth of a finitely generated $VIC(\mathbb{Z})$-module is at always bounded by some exponential function $C^n$.  The case where the dimensions grow even slower is of particular importance and various classes of such modules have appeared in the literature. In this section we will give a characterization of $VIC(\mathbb{Z})$-modules of slow dimension growth, showing that several different conditions are equivalent.

Before stating the characterization let us recall the notion of finite polynomial degree for $VIC(R)$-modules due to Randal-Williams and Wahl (\cite{RWW} Definition 4.10). We say that a $VIC(R)$-module $V$ has polynomial degree $-\infty$ in degrees $>N$ if  $V_n = 0$ for all $n>N$. Then inductively for $d \ge 0$ we say that $V$ has polynomial degree $\le d$ in degrees $>N$ if

\begin{enumerate}

\item The kernel of the natural map $V \to SV$ is zero in degrees $>N$.

\item The cokernel of this map $V \to SV$ has polynomial degree $\le d-1$  in degrees $>N$.

\end{enumerate}

We will say that $V$ has finite polynomial degree in sufficiently large degree if there exists a $d$ and an $N$ such that $V$ has polynomial degree $\le d$ in degrees $>N$.  The following proposition characterizes $VIC(\mathbb{Z})$-modules of small dimension growth.

\begin{proposition}

Let $V$ be a finitely generated $VIC(\mathbb{Z})$-module over $\mathbb{C}$.  The following conditions are equivalent:

\begin{enumerate}[(1)]

\item $\dim(V_n) < 2^n-2$ for all $n \gg0$.

\item There exists a polynomial $p(x) \in \mathbb{Q}[x]$ such that $\dim(V_n) = p(n)$ for all $n \gg 0$.

\item $V$ is pointwise algebraic in sufficiently large degrees.

\item $V$ is finitely generated as an $FI$-module.

\item $V$ is has finite polynomial degree in sufficiently large degree.

\end{enumerate}

\end{proposition}

\noindent \textbf{Proof:}  $(1) \Rightarrow (3)$ is immediate from Corollary \ref{lowdimalg}. Next, Corollary \ref{algextend} tells us that a pointwise algebraic $VIC(\mathbb{Z})$-module extends to a pointwise algebraic $VIC(\mathbb{C})$-module in large enough degrees. Patzt's theorem on pointwise algebraic $VIC(\mathbb{C})$-modules tells us that any such $VIC(\mathbb{C})$-module $V$ is multiplicity stable meaning that $$V_n \cong \bigoplus_{\lambda^+,\lambda^-} V(\lambda^+,\lambda^-)^{c_{\lambda^+,\lambda^-}}$$ where the multiplicities ${c_{\lambda^+,\lambda^-}} $are independent of $n$ for $n \gg 0$ and nonzero for only finitely many ${(\lambda^+,\lambda^-)}$. Therefore $(3)\Rightarrow(2)$ since for any bipartition $(\lambda^+,\lambda^-)$ the dimension of $V(\lambda^+,\lambda^-)$ is given by a polynomial in $n$ of degree $|\lambda^+| + |\lambda^-|$ by the Weyl dimension formula.  Of course $(2) \Rightarrow (1)$ is obvious, so $(1)$, $(2)$, and $(3)$ are equivalent.

$(5) \Rightarrow (2)$ is clear from the definition of polynomial degree, but under these conditions we claim the converse is also true. Since $V$ is assumed to be finitely generated Corollary \ref{inject} tells us that the kernel of the map $V \to SV$ is zero in sufficiently large degree, so $V$ has finite polynomial degree in sufficiently large degree if and only if the cokernel of $V \to SV$ has finite polynomial degree. So now by induction on polynomial degree we see that the class of $VIC(\mathbb{Z})$-modules polynomial degree $\le d$ in sufficiently large degree is the same as the class of $VIC(\mathbb{Z})$-modules such that $\dim(V_n)$ is given by a polynomial of degree at most $d$
. Finally, the equivalence of $(4)$ and $(5)$ is Theorem 3.25 of $\cite{MPW}$ . $\square$

\subsection{$VIC(\mathbb{Z})$-modules outside characteristic zero}

We'd like to say a bit about $VIC(\mathbb{Z})$-modules over other rings and in particular in positive characteristic. It seems reasonable to conjecture that Noetherianity holds for finitely generated $VIC(\mathbb{Z})$-modules that are pointwise finitely generated over a Noetherian ring, but this seems out of reach at the moment.  

In general we don't have a version superrigidity over arbitrary ground rings (although some things can be deduced via base change). We will note however that the notion of finite type representations (those that factor through a finite quotient) and their depth make sense over any commutative ground ring. Moreover, the proof of Theorem \ref{fintypevic} is group theoretic and also works over any commutative ground ring. As a corollary we have the following Noetherianity result for pointwise finite type $VIC(\mathbb{Z})$-modules: 

\begin{corollary}

If $V$ is a finitely generated $VIC(\mathbb{Z})$-module over a Noetherian ring $R$ which is pointwise finitely generated over $R$ and of finite type then any $VIC(\mathbb{Z})$-submodule $W \subseteq V$ is finitely generated. 

\end{corollary}

\noindent \textbf{Proof:} Theorem \ref{fintypevic} tells us that any such module is actually a $VIC^\pm(\mathbb{Z}/\ell \mathbb{Z})$-module for some $\ell$ and Putman and Sam have shown $VIC^\pm(\mathbb{Z}/\ell \mathbb{Z})$-modules over a Noetherian ring are Noetherian (\cite{PS} Theorem D). $\square$

\medskip

In some situations we care about this condition of being pointwise finite type will be automatic. In particular if we knew that every representation of $SL_n(\mathbb{Z})$ that is finitely generated over a ring $R$ for $n >2$ was of finite type this would apply automatically.  We'll call such rings \emph{SL-finite}.

\begin{corollary}

If $V$ be a finitely generated $VIC(\mathbb{Z})$-module which is pointwise finitely generated over a $SL$-finite Noetherian ring $R$ then any $VIC(\mathbb{Z})$-submodule is finitely generated.

\end{corollary}

This condition of being $SL$-finite may seem a bit strange but we'll note that it clearly includes all finite rings. The following lemma will tell us it also includes all fields of characteristic $p > 0$ (and therefore $VIC(\mathbb{Z})$-modules  over a field of characteristic $p$ are Noetherian). It may also be of independent interest.

\begin{lemma}
Let $k$ be a field of characteristic $p > 0$.  For $n>2$ any finite dimensional representation $SL_n(\mathbb{Z}) \to GL_N(k)$ factors through $SL_n(\mathbb{Z} / \ell \mathbb{Z})$ for some $\ell$.

\end{lemma}

\noindent \textbf{Proof:} First let's assume without loss of generality that $k$ is algebraically closed. If $k = \bar{\mathbb{F}}_p$ then since $SL_n(\mathbb{Z})$ is finitely generated its image must be contained in $GL_N(\mathbb{F}_{p^m})$ for some $m$, but of course this is a finite group so the congruence subgroup property tells us the map must factor through $SL_n(\mathbb{Z} / \ell \mathbb{Z})$ for some $\ell$. We will reduce the general case to this via a bit of model theory and some dimension bounds.

 Since $SL_n(\mathbb{Z})$ is finitely presented the condition of there existing an $N$ dimensional representation of $SL_n(\mathbb{Z})$ where $E$ acts with order $> \ell$ is a first order sentence in the language of algebraically closed fields. We know that all algebraically closed fields of characteristic $p$ are elementary equivalent, so this is true over $k$ if and only if it is true over $\bar{\mathbb{F}}_p$.  So if there is an $N$ dimensional representation over $k$ which is not of finite type then $E$ has infinite order and therefore there must be finite type $N$ dimensional representations of arbitrarily high depth over $\bar{\mathbb{F}}_p$.

So now it is enough to show that for fixed $n > 2$ and $N$ there are finitely many values of $\ell$ such that there is a finite type depth $\ell$ representation of $SL_n(\mathbb{Z})$ of dimension at most $N$ over $k$. By the Chinese remainder theorem it suffices to do this when $\ell = q^m$ is a prime power.

 If $q \ne p$ this follows from the bounds of Bardestani and Mallahi-Karai stated in Theorem \ref{pkbound}.  They state these bounds over $\mathbb{C}$ but their proof only needs that the representations of a certain subgroup isomorphic to $(\mathbb{Z} / q^m \mathbb{Z})^{n-1}$ are direct sums of one dimensional characters, which holds over any algebraically closed field of characteristic prime to $q$.

If $q=p$ this argument breaks down and indeed $SL_n(\mathbb{Z}/p\mathbb{Z})$ has representations in characteristic $p$ that are much smaller than the bounds in characteristic zero (for example the obvious $n$ dimensional one). For our purposes though we just need to give any bound. We know that $E$ must act with order $p^m$, and the only elements of $GL_N(k)$ with order a power of $p$ are unipotent elements, the largest such order of which is when $m = \lceil \text{log}_p(N) \rceil$. $\square$

\section{Questions and future directions of work} \label{future}

\subsection{Variations}

This paper is in part a proof-of-concept both that rigidity results can be applied to control families of representations and that the type of finiteness and stability conditions arising in representation stability can be used to extend superrigidity to infinite rank settings. We'd like to suggest some further directions along that vein.

\subsubsection*{Symplectic groups}

Instead of sequences of $GL_n(\mathbb{Z})$-representations we could look at sequences of $Sp_{2n}(\mathbb{Z})$-representations. Many of the ingredients used in this paper also hold in this case:

\begin{itemize}

\item These groups also satisfy the congruence subgroup property.

\item Superrigidity again gives that irreducibles are tensor products of algebraic and finite type parts.

 \item There are similar dimension bounds for representations of $Sp_{2n}(\mathbb{Z} / p^k \mathbb{Z})$ (again found in \cite{TZ} and \cite{BMK}).
 
 \item There category $SI(\mathbb{Z})$ is a natural symplectic version of $VIC(\mathbb{Z})$.
 
\item  Patzt has also proven a representation stability result for pointwise algebraic $SI(\mathbb{C})$-modules via a similar analysis of the branching rules.

\end{itemize}
It's likely that all of the main results here will go through mutatis mutandis in this setting.

\subsubsection*{Other rings of integers}

People have also looked at $VIC(R)$-modules for other rings $R$. One particularly important case is when $R = \mathcal{O}$ is the ring of integers in some number field. It would be interesting to generalize the results of this paper to that case, to understand the structure of $VIC(\mathcal{O})$-modules and to prove a superrigidity result for $GL_{\infty}(\mathcal{O})$.

 The classification of irreducible representations via superrigidity for $GL_n(\mathcal{O})$ is more complicated than for $GL_n(\mathbb{Z})$, for one thing we should expect ``algebraic" components coming from all embeddings of $\mathcal{O}$ into $\mathbb{C}$, and the structure of the profinite completion of $GL_n(\mathcal{O})$ will depend on the number theoretic properties of $\mathcal{O}$ more subtly than for $\mathbb{Z}$. 

We'll also note that some of the nice properties of $VIC(R)$-modules are only known when $R$ is a PID, which presents another technical hurdle.  We suspect though that, especially with the additional assumption of pointwise finite dimensionality, $VIC(\mathcal{O})$-modules will have many of these properties as well, even if $\mathcal{O}$ is not a PID.

\subsubsection*{$SO(p,q)$ and infinite dimensional hyperbolic manifolds}

Another interesting direction to look in is at families of lattices in special orthogonal groups $SO(p,q)$ where we fix either $p$ or $q$ and let the other grow. In particular if we fix $p = 1$ it would be especially interesting to consider compatible families of cocompact lattices $\Gamma_q \subset SO(1,q)$ and to study the corresponding infinite dimensional hyperbolic manifolds $\mathbb{H}^\infty / \Gamma_\infty$. This is rapidly leaving the realm of things the author is qualified to speculate on though, so we won't say anything more specific.

\subsection{Possible improvements}

This paper explores pointwise finite dimensional $VIC(\mathbb{Z})$-modules and sheds light on to certain aspects some of them, but it is definitely not the complete picture. Here are just a few possible places where one could improve our understanding. 

\subsubsection*{Degree bounds}

Our main interest in this paper was on understanding the eventual behavior of finitely generated pointwise finite dimensional $VIC(\mathbb{Z})$-modules.  For many applications though it is useful to know not just what the eventual behavior is, but to understand at what degrees the stable behavior actually starts.  Especially important is bounding the injectivity and surjectivity degrees in terms of the presentation degree.

For $VIC(\mathbb{F}_q)$-modules Miller and Wilson have obtained some very tight bounds for these sort of questions (in \cite{MW}), it is likely that some of their analysis can be carried over to this setting as well.

\subsubsection*{General Noetherianity}

We now have a Noetherianity statement for $VIC(\mathbb{Z})$-modules over $\mathbb{C}$ via superrigidity and semisimplicity, as well as one over finite rings and fields of characteristic $p$ using the $SL$-finiteness property.  It seems reasonable that Noetherianity should also hold in greater generality.  Explicitly, we expect that the following holds:

\begin{conjecture}

If $V$ is a finitely generated $VIC(\mathbb{Z})$-module which is pointwise finitely generated over a Noetherian ring $R$ then any $VIC(\mathbb{Z})$-submodule is also finitely generated.

\end{conjecture}

\subsubsection*{Propagating finite dimensionality}

We have been assuming pointwise finite-dimensionality as an external condition on our $VIC$-modules. It would be interesting to look at conditions which would guarantee this, or to see if one can propagate finite dimensionality from low degrees onward. In particular we'd like to pose the following question:

Suppose $V$ is a $VIC(\mathbb{Z})$-module which is generated in degrees $\le d$, and suppose for some $n > d$ ($n = d+1$ is especially interesting) we know that $V_n$ is finite dimensional.  Is it true that $V_m$ is finite dimensional for all $m \ge n$? If not, what if instead we know that two adjacent terms $V_n$ and $V_{n+1}$ are finite dimensional? How about if $n$ is sufficiently large compared to $d$?

We will note that if we additionally assume that the representations are pointwise finite type then we can do this (over any ground ring): If $V$ is generated in degree $\le d$ and is finite dimensional and of finite type with depth dividing $\ell$ in degrees $d$ and $d+1$ then Lemma \ref{depthprop} implies from that point on $V$ is a $VIC(\mathbb{Z} /\ell \mathbb{Z})$-module, and for those of course finite generation implies pointwise finite dimensionality.

\subsubsection*{Decoupling the different local pieces}

Theorem \ref{main} can be thought of as decoupling the algebraic part from the finite type part of a pointwise finite $VIC(\mathbb{Z})$-module.  It seems likely that one could similarly decouple the contributions from the factors coming from different primes. 

 In particular, one way to obtain a finitely generated pointwise finite dimensional and finite type $VIC(\mathbb{Z})$-module is to take some finitely generated pointwise smooth $VIC(\mathbb{Z}_p)$-modules for different primes $p$ and then tensor them together.  Is it true that every finitely generated pointwise finite dimensional and finite type $VIC(\mathbb{Z})$-module has a filtration by modules obtained this way?

 \subsubsection*{Multiplicity stability for $VIC(\mathbb{Z}/p^k\mathbb{Z})$-modules}

Gan and Watterlond \cite{GW} have shown that $VIC(\mathbb{F}_p)$-modules over $\mathbb{C}$ exhibit representation stability in the sense that the multiplicities stabilize with respect to certain natural labelings of the representations of $GL_n(\mathbb{F}_p)$. The irreducible components of a $VIC(\mathbb{Z})$-module of depth not dividing $p$ form a $VIC$-submodule, so looking at the quotient by this, the Gan and Watterlond result combined with Theorem \ref{main} gives us multiplicity stability for the depth $p$ components of a pointwise finite dimensional $VIC(\mathbb{Z})$-module.

It seems reasonable to expect that something analogous should hold for $VIC(\mathbb{Z}/p^k\mathbb{Z})$-modules (and hence for pointwise smooth $VIC(\mathbb{Z}_p)$-modules), however we are unaware of any sort of nice labeling of the representations of $GL_n(\mathbb{Z}/p^k\mathbb{Z})$ to work with. Explicitly we'll give the following conjecture:

\begin{conjecture}

For every $p, k$ there exists a compatible collection of labelings of the irreducible complex representations of $GL_n(\mathbb{Z}/p^k\mathbb{Z})$ for all $n$ such that finitely generated $VIC(\mathbb{Z}/p^k\mathbb{Z}))$-modules exhibit multiplicity stability with respect to this labeling. 

\end{conjecture}

 If this holds, and if the different local contributions can be decoupled as suggested above, one might may be able to formulate a multiplicity stability result for pointwise finite dimensional $VIC(\mathbb{Z})$-modules.

\subsection{The inverse-transpose automorphism}

The following line of thought didn't lead anywhere useful for this project, but it was somewhat interesting to work out so we've included a brief summary.

\medskip

Recall that $GL_n$ has an outer automorphism given by sending a matrix $A$ to its inverse-transpose $A^{-1 \top}$.  It turns out this actually extends to an automorphism of the entire category $VIC$, and in some ways makes it seem more natural in this light.

Recall that we can think of $VIC(R)$ as being the category of free finite rank $R$-modules with split injections.  That is, morphisms from $M$ to $N$ are pairs $(f,g)$ where $f\colon M \to n$ is an injective $R$-linear map, and $g\colon N \to M$ is a an $R$-linear map such that $g \circ f = \text{Id}_M$.  This map $g$ is the same as the data of a complementary subspace $M^\perp$ such that $N = Im_f(M) \oplus M^\perp$, by letting $M^\perp = \text{ker}(g)$. 

The automorphism of $VIC$ sends $M$ to its dual $M^*$, and a split injection $(f,g)$ from $M$ to $N$ to $(g^*,f^*)$ from $M^*$ to $N^*$.  In particular it sends an endomorphism $(f,f^{-1})$ of $M$ to the endomorphism $(f^{-1*}, f^*)$ of $M^*$, which in terms of a basis for $M$ and the dual basis for $M^*$ sends a matrix $A$ to $A^{-1 \top}$.

We can twist a $VIC$-module $V$ by precomposing with this automorphism to obtain a new $VIC$-module $\tilde{V}$.  We can think of $V$ and $\tilde{V}$ as consisting of the same vector spaces with the same underlying set of linear maps between them, but with the labeling of which map in $VIC$ corresponds to which linear map twisted by this automorphism.  In particular, it's clear that $\tilde{V}$ is finitely generated if and only if $V$ is.

So what does this do to representations?  Well for algebraic representations the answer is simple:  If I precompose an algebraic representation of $GL_n(\mathbb{C})$ with the inverse-transpose map, the new representation is isomorphic to the dual of the original\footnote{If one takes a representation $G \to GL_n(\mathbb{C})$ and then \emph{post}composes it with the inverse-transpose map, then the resulting representation is always dual to the original representation of $G$. This is well known, but is not what is happening here.}. This follows from the fact that every matrix over $\mathbb{C}$ (or any field) is conjugate to its own transpose, so the character of an element $g$ on the twisted representation is the character of $g^{-1}$ on the original representation and of course the dual representation has this same property. Since algebraic representations in characteristic zero are semisimple, having the same character implies these representations are isomorphic.

 In particular, if $V$ is an algebraic $VIC(\mathbb{C})$-module then $\tilde{V}$ is pointwise dual to $V$ meaning $\tilde{V}_n \cong V_n^*$.  Note though, these isomorphisms are very non-canonical and the maps between the $\tilde{V}_n$'s are not related in an obvious way to the maps between the $V_n$'s.

What about for non-algebraic representations?  Well it turns out that for $n \ge 3$ not every matrix $GL_n(\mathbb{Z})$ is conjugate to its transpose, and in fact there are matrices over $\mathbb{Z}/p^2\mathbb{Z}$ which are not conjugate to their transposes (see \cite{PSS1}). In particular this implies there exists representations of $GL_n(\mathbb{Z})$ for which the dual and the inverse-transpose twist representations are non-isomorphic, as their characters will be different on these non self-transpose conjugacy classes.

We'll also note that the isomorphism for algebraic representations also breaks down in positive characteristic, but for different reasons. Irreducible algebraic representations of $GL_n(\bar{\mathbb{F}}_p)$ are still determined by their characters, so the same reasoning tells us that the twist is isomorphic to the dual. However if say we have a non-trivial extension of two irreducibles with $V_1$ as a subrepresentation and $V_2$ as a quotient. The twist has $V_1^*$ as a subrepresentation and $V^*_2$ as a quotient, whereas the dual has $V_2^*$ as a subrepresentation and $V^*_1$ as a quotient.

\subsection{A conjecture on representation growth}

Finally we'd like to close out with a conjecture about the number of irreducible finite type representations $SL_n(\mathbb{Z})$ of small size. For fixed $n$ there are bounds due to Aizenbud and Avni \cite{AA} on the number of irreducible representations of $SL_n(\mathbb{Z}_p)$ of dimension at most $d$, which depend on $d$ but not $n$ or $p$.  

However we'll note that if we fix $d$ then for $n>>0$ there are no nontrivial representations $SL_n(\mathbb{Z}_p)$ with dimension less than $d$.  Similarly, if we fix $d$ and $n$ then Theorem \ref{pkbound} says that any representation of $SL_n(\mathbb{Z}_p)$ with dimension less than $d$ must have sufficiently small depth (and hence there are only finitely many such representations over all $p$ combined). 

So while Aizenbud-Avni bounds are great for controlling the number irreducible representations of $SL_n(\mathbb{Z}_p)$ for fixed $n$ as the dimension gets large, they don't give much insight into the counting low dimensional representations.  We'd now like to conjecture a sort of logarithmic version of the Aizenbud-Avni counts for the low dimensional representations:

\begin{conjecture}

For every constant $C > 0$ there exists a bound $b(C)$ such that for all $n, p$  there are at most $b(C)$ irreducible smooth representations of $SL_n(\mathbb{Z}_p)$ having dimension less than $C^n$. 

\end{conjecture}

Of course if this conjecture is true it would be also very interesting to understand the rate of growth of $b(C)$ as a function of $C$, or to package this information into a sort of ``stable" representation zeta function and understand its convergence properties as Aizenbud and Avni do for fixed $n$, but we won't venture any specific conjectures in that direction at this time.

\bibliographystyle{plain}
\bibliography{bib}

\end{document}